\documentclass[onefignum,onetabnum]{siamart220329} 

\usepackage{lipsum}
\usepackage{amsfonts}
\usepackage{graphicx}
\usepackage{epstopdf}
\usepackage{algorithmic}
\ifpdf
  \DeclareGraphicsExtensions{.eps,.pdf,.png,.jpg}
\else
  \DeclareGraphicsExtensions{.eps}
\fi

\newsiamremark{remark}{Remark}
\newsiamremark{hypothesis}{Hypothesis}
\crefname{hypothesis}{Hypothesis}{Hypotheses}
\newsiamthm{claim}{Claim}

\newcommand{\R}{\mathbb{R}}

\newcommand{\matrx}[1]{\mathbf{#1}}
\renewcommand{\vec}[1]{\mathbf{#1}}
\newcommand{\A}{\matrx{A}} 
\newcommand{\Sm}{\matrx{S}} 
\newcommand{\inex}{\mathrm{in}}
\newcommand{\ex}{\mathrm{ex}}
\newcommand{\prev}{\mathrm{prev}}
\newcommand{\new}{\mathrm{new}}

\newcommand{\rv}[1]{\textcolor{black}{#1}} 
\newcommand{\rev}[1]{\textcolor{black}{#1}}

\usepackage{tikz}
\usepackage{pgfplots}
\usetikzlibrary{positioning}
\usetikzlibrary{plotmarks}
\newlength\figureheight
\newlength\figurewidth

\usepackage{amssymb}

\headers{Approximate coarsest-level solves in multigrid}{P. Vacek, E. Carson, and K. M. Soodhalter}

\title{The effect of approximate coarsest-level solves on the convergence of multigrid V-cycle methods\thanks{
 We acknowledge funding from Charles University PRIMUS project no. PRIMUS/19/SCI/11, the grant SVV-2023-260711, Charles University Research Centre program No. UNCE/24/SCI/005, the Exascale Computing Project (17-SC-20-SC), a collaborative effort of the U.S. Department of Energy Office of Science and the National Nuclear Security Administration, and by the European Union (ERC, inEXASCALE, 101075632). Views and opinions expressed are those of the
authors only and do not necessarily reflect those of the European Union or the European Research Council. Neither the European Union nor the granting authority can be held responsible for them.
}}

\author{Petr Vacek\thanks{Department of Numerical Mathematics, Faculty of Mathematics and Physics, Charles University, (\email{\{vacek,carson\}@karlin.mff.cuni.cz}).}
\and Erin Carson\footnotemark[2]
\and Kirk M. Soodhalter\thanks{School of Mathematics, Trinity College Dublin, (\email{ksoodha@maths.tcd.ie}).}}

\usepackage{amsopn}

\ifpdf
\hypersetup{
  pdftitle={The effect of approximate coarsest-level solves on the convergence of multigrid V-cycle methods},
  pdfauthor={P. Vacek, E. Carson, and K. M. Soodhalter}
}
\fi

\begin{document}

\maketitle

\begin{abstract}
The multigrid V-cycle method is a popular method for solving systems of linear equations. It computes an approximate solution by using smoothing on fine levels and solving a system of linear equations on the coarsest level. 
Solving on the coarsest level depends on the size and difficulty of the problem.
If the size permits, it is typical to use a direct method based on LU or Cholesky decomposition. In settings with large coarsest-level problems, approximate solvers such as iterative Krylov subspace methods, or direct methods based on low-rank approximation, are often used. The accuracy of the coarsest-level solver is typically determined based on the experience of the users with the concrete problems and methods.

\rv{
In this paper we present an approach to analyzing the effects of approximate coarsest-level solves on the convergence of the V-cycle method for symmetric positive definite problems.
Using these results, we derive coarsest-level stopping criterion through which we may control the
difference between the approximation computed by a V-cycle method with approximate coarsest-level solver and the approximation which would be computed if the coarsest-level problems were solved exactly.
The coarsest-level stopping criterion may thus be set up such that the V-cycle method converges to a chosen finest-level accuracy in (nearly) the same number of V-cycle iterations as the V-cycle method with exact coarsest-level solver.
We also utilize the theoretical results to discuss how the convergence of the V-cycle method may be affected by the choice of a tolerance in a coarsest-level stopping criterion based on the relative residual norm.
}
\end{abstract}



\section{Introduction}
Multigrid methods \cite{Brandt2011,Briggs2000,Trottenberg2001,Hackbusch2016} are frequently used when solving systems of linear equations, and can be applied either as standalone solvers or as preconditioners for iterative methods. There are two types of multigrid; \emph{geometric}: wherein the hierarchy of systems is obtained by discretizations of an infinite dimensional problem on a sequence of nested meshes; and \emph{algebraic}: wherein the coarse systems are assembled based on the algebraic properties of the matrix. Within each multigrid cycle, the approximation is computed using smoothing on fine levels and solving a system of linear equations on the coarsest level. Smoothing on the fine levels is typically done via a few iterations of a stationary iterative method. The particular solver used for the problem on the coarsest level depends on its size and difficulty. If the size of the problem permits, it is typical to use a direct solver based on LU or Cholesky decomposition.

In this text, we focus on settings where the problem on the coarsest level is large and the use of direct solvers based on LU or Cholesky decomposition may be ineffective or impossible to realize.
Such settings may arise, for example, when using geometric multigrid methods to solve problems on complicated domains. The mesh associated with the coarsest level must resolve the domain with certain accuracy. This can yield a large number of degrees of freedom. One possible solution to this issue is to solve the coarsest-level problem using algebraic multigrid, which can introduce additional coarse levels that are not related to the geometry of the problem.

Another setting where large coarsest-level problems may be present is when we use multigrid methods on parallel computers. In parallel computing, the degrees of freedom are assigned to different processors or accelerators. The computation is done in parallel on the individual processors and the results are communicated between them. A challenge for effective parallel implementation of multigrid methods is that the amount of computation on coarse levels decreases at a faster rate than the amount of communication; see e.g., the discussion in the introduction of \cite{Buttari2022}. \rv{One possible solution is to treat this issue by redistribution of the coarse-level problems to a smaller number of processors; see e.g., \cite{Gahvari2013,May2016,Reisner2018}.
Another solution may be to use communication-avoiding methods on the coarse levels; see e.g., \cite{Williams2014}.}

In this paper, we instead consider treating the still large-scale coarsest-level problem by solving \emph{inexactly}. Frequently used solvers for large scale coarsest-level problems include Krylov subspace methods and direct approximate solvers; see, e.g., \cite{Huber2019}, where the author considers the preconditioned conjugate gradient method, or \cite{Buttari2022}, where the authors study the use of a block low-rank (BLR) low precision direct solver.
These solvers approximate the coarsest-level solution to an accuracy which is determined by the choice of a stopping criteria or affected by the choice of the low-rank threshold and finite precision. These parameters are often chosen in practice based on the experience of the user with concrete problems and methods with the goal of balancing the cost of the coarsest-level solve and the total number of V-cycles required for convergence. In \Cref{section.motivating-experiment} we present a motivating numerical experiments, which illustrate how the choice of the accuracy of the coarsest-level solver may affect the convergence of the multigrid V-cycle method.

A general analysis of the effects of the accuracy of the coarsest-level solver on the convergence behaviour of multilevel methods is, to our knowledge, not present in the literature. Multigrid methods are typically analyzed under the assumption that the problem on the coarsest level is solved exactly; see, e.g., \cite{Yserentant1993,Xu1992}.
An algebraic analysis of perturbed two grids methods and its application to the analysis of other multigrid schemes with approximate coarsest-level solvers can be found in \cite{Notay2007,Xu2022}. \rv{The authors derive estimates of the worst-case convergence rate of the methods. The results are, however, obtained under the assumption} that the action of the solver on the coarsest level can be expressed using a symmetric positive definite matrix. This is not true for frequently used solvers, e.g., for a Krylov subspace method stopped using a relative residual stopping criterion.
A more general setting is considered in the paper \cite{McCormick2021}, which presents the first analysis of mixed precision multigrid solvers. The authors assume that the action of the solver on the coarsest level can be expressed using a non-singular matrix.

In this paper, we propose an approach to algebraically analyze the effect of approximate coarsest-level solves in the multigrid V-cycle method for symmetric positive definite (SPD) problems. 
The main methodology of our approach is to view the inexact V-cycle (inV-cycle) method as a perturbation of the exact V-cycle (exV-cycle) method in the following sense. We express the error of the approximation computed by one V-cycle with an approximate coarsest-level solver as the error of the approximation computed by one V-cycle with an exact coarsest-level solver plus the difference of the two approximations. We show that the difference can be expressed as a matrix times the error of the coarsest-level solver. The matrix describes how the error from the coarsest level is propagated to the finest level. \rv{Moreover, we consider two assumptions on the accuracy of the coarsest-level solver:
a \emph{relative} assumption, where the error of the coarsest-level solver is less than a factor of the error of the previous finest-level approximation, and an \emph{absolute} assumption, where the error of the coarsest-level solver is less than a certain constant.
Based on the relative assumption we derive an estimate on the convergence rate of the inV-cycle method and discuss its uniform convergence.
Utilizing the absolute assumption we get an estimate on the difference between the approximation computed by the inV-cycle method and the exV-cycle method after a number of V-cycle iterations.
}
The analysis is done assuming exact arithmetic computations, aside from the computation of the coarsest level solutions. The model is agnostic about what coarsest-level solver is used; we only assume that the error on the coarsest level satisfies certain assumptions.

The paper is organized as follows. In \Cref{sec:Vcycle} we establish the notation, state the V-cycle method and present a motivating numerical experiments, which illustrate that the choice of the accuracy of the coarsest-level solver can significantly affect the convergence of the V-cycle method.
In \Cref{sec:inV-cycle_analysis} we present an analysis of the V-cycle method with an approximate coarsest-level solver. 
The results are applied to describe the possible effects of the choice of the tolerance in a coarsest-level relative residual stopping criterion in \Cref{sec:effects_rel_res}. 
New stopping criteria based on the absolute coarsest-level accuracy assumption are derived in \Cref{sec:new_stopping_criteria}. 
Finally, we present a series of numerical experiments illustrating the obtained results in \Cref{sec:numerical_experiments}. The text closes with conclusions and discussion of open problems in \Cref{sec:conclusion}.

\section{Notation and motivating experiments}
\label{sec:Vcycle}
We study the multigrid V-cycle method for finding an approximate solution of the following problem. Given an SPD matrix $\matrx{A}\in \R^{n \times n}$ and a right-hand side vector $\vec{b} \in \R^{n}$ find the vector $\vec{x} \in \R^{n}$ such that
\begin{equation*}
\matrx{A}\vec{x} = \vec{b}.
\end{equation*}
We consider a hierarchy of $J+1$ levels numbered from zero to $J$\rv{, where level zero is the coarsest level and level $J$ the finest level}. Each level contains a system matrix $\matrx{A}_j\in\R^{n_j \times n_{j}}$, with $\matrx{A}_J = \A$. Information is transferred between the $(j-1)$th level and the $j$th level using a full rank prolongation matrix $\matrx{P}_{j}\in\R^{n_j \times n_{j-1}}$, respectively its transpose. We assume that the system matrices and the prolongation matrices satisfy the so called \emph{Galerkin condition}, i.e.,
\begin{equation}\label{eq:galerkin_condition}
\matrx{A}_{j-1} =  \matrx{P}_{j}^{\top} \matrx{A}_{j} \matrx{P}_{j}, \quad j=1\ldots,J.
\end{equation}
We use the notation $ \matrx{A} _{0:j}$, for the sequence of matrices $\matrx{A}_{0}, \ldots ,\matrx{A}_{j}$.
Let $\| \cdot\|$ denote the Euclidean vector norm and  let $\| \cdot \|_{\matrx{A}_j} = \| \matrx{A}^{\frac{1}{2}}_j \cdot \|$ denote the $\matrx{A}_j$ vector norm, also called the energy norm. We use the same notation for the matrix norms generated by the associated vector norms.
Let $\matrx{I}_j\in\R^{n_j \times n_{j}}$ denote the identity matrix on the $j$th level.

We assume that the pre- and post- smoothing on levels $j=1,\ldots,J$ can be expressed in the form 
\begin{equation*}
    \vec{v}_j  = \vec{v}_j +  \matrx{M}_j (\vec{f}_j - \matrx{A}_j \vec{v}_j) \quad \text{and}
    \quad \vec{v}_j  = \vec{v}_j +  \matrx{N}_j (\vec{f}_j  - \matrx{A}_j \vec{v}_j),
\end{equation*}
respectively, where $\vec{v}_j$ and $\vec{f}_j$ are an approximation and a right-hand side on the $j$th level and $\matrx{M}_j \in \R^{n_j \times n_{j}}$ and $\matrx{N}_j\in \R^{n_j \times n_{j}}$ are non-singular matrices satisfying 
\begin{equation}\label{eq:smoothing_convergence_Anorms}
\| \matrx{I}_j - \matrx{M}_j\matrx{A}_j \|_{\matrx{A}_j}<1
\quad \text{and} \quad
\| \matrx{I}_j - \matrx{N}_j\matrx{A}_j \|_{\matrx{A}_j}<1.
\end{equation}
This assumption yields monotone convergence of the smoothers as standalone solvers in the $\A_j$-norms. 
Frequently used smoothers, e.g., a few iterations of a classic stationary iterative method such as damped Jacobi or Gauss-Seidel,  typically satisfy these assumptions; see, e.g., the discussion in \cite[p.~293]{Yserentant1993} or \cite{Xu1992}.
We also consider multilevel schemes, where either pre- or post- smoothing is not used,  i.e., where formally either $\matrx{M}_j$, $j=1,\ldots,J$ or $\matrx{N}_j$, $j=1,\ldots,J$, are zero matrices.

Given an approximation $\vec{x}^{\prev}$ to the solution $\vec{x}$, the approximation after one iteration of the V-cycle method is computed by calling \Cref{alg:V-cycle} as  (see, e.g., \cite[pp.~47--48]{Trottenberg2001})
\begin{equation*}
\vec{x}^{\new} = \mathbf{V}( \matrx{A} _{0:J}, \matrx{M}_{1:J}, \matrx{N}_{1:J}, \matrx{P}_{1:J}, \vec{b},\vec{x}^{\prev},J).
\end{equation*}
We distinguish between the \emph{exV-cycle method} and the \emph{inV-cycle method} based on whether the coarsest-level problem is solved exactly or not. 
\begin{algorithm}
\caption{V-cycle scheme, $\mathbf{V}( \matrx{A} _{0:j}, \matrx{M}_{1:j}, \matrx{N}_{1:j}, \matrx{P}_{1:j}, \vec{f}_{j}$, $\vec{v}^{[0]}_{j},j)$.}
\label{alg:V-cycle}
\rv{
\begin{algorithmic}
\IF{$ j \neq 0 $}
\STATE {$\vec{v}^{[1]}_{j} = \vec{v}^{[0]}_j + \matrx{M}_j (\vec{f}_{j} - \matrx{A}_{j} \vec{v}^{[0]}_j)$} \COMMENT{pre-smoothing}
\STATE{$\vec{f}_{j-1} = \matrx{P}_{j}^{\top} (\vec{f}_{j-1} - \matrx{A}_{j} \vec{v}^{[1]}_j)$ } \COMMENT{restriction}
\STATE{$\vec{v}^{[2]}_{j-1} = \mathbf{V}(\matrx{A}_{0:j-1}, \matrx{M}_{1:j-1}, \matrx{N}_{1:j-1}, \matrx{P}_{1:j-1}, \vec{f}_{j-1}$, $\vec{0}$, $j-1)$}
\STATE{$\vec{v}^{[3]}_{j} = \vec{v}^{{[1]}}_{j} + \matrx{P}_{j}  \vec{v}^{[2]}_{{j-1}}$ }\COMMENT{coarse grid correction}
\STATE{ $\vec{v}^{[4]}_{j}=\vec{v}^{{[3]}}_{j} +  \matrx{N}_j(\vec{f}_j  - \matrx{A}_j\vec{v}^{[3]}_{{j}})$} \COMMENT{post-smoothing}
\RETURN{$\vec{v}^{[4]}_{j}$}
\ELSE
\RETURN{(approximate) solution of the problem $\matrx{A}_0 \vec{v}_{0} =  \vec{f}_0$}
\ENDIF
\end{algorithmic}
}
\end{algorithm}

\subsection{Motivating experiments}\label{section.motivating-experiment}
We illustrate the relevance of the forthcoming analysis with numerical experiments, which demonstrate how the choice of the accuracy of the coarsest-level solve affects the convergence of the V-cycle method.

\rv{
We consider a second order elliptic PDE of the form
\begin{equation*}
    - \nabla \cdot (  k(x) \nabla u ) = f \quad \text{in } \Omega, \qquad 
    u = 0   \quad \text{on } \partial \Omega,
\end{equation*}
where $f \equiv 1$ and $\Omega = (0,1)\times(0,1)$. We consider two variants of the problem based on the coefficient function $k:\Omega \rightarrow \mathbb{R}$, ``Poisson'' with $ k\equiv 1$ and ``jump-1024'' with 
\begin{equation*}
    k(x) = 
    \begin{cases}
    1024, \quad x \in \left(0,\frac{1}{2} \right) \times \left(0,\frac{1}{2} \right) \cup \left( \frac{1}{2},1 \right)\times \left( \frac{1}{2},1\right), \\
    1, \quad x \in \left(0,\frac{1}{2}\right)\times\left(\frac{1}{2},1\right) \cup \left(\frac{1}{2},1\right)\times\left(0,\frac{1}{2}\right).
    \end{cases}
\end{equation*}
}

The problems are discretized using the Galerkin finite element (FE) method with continuous piecewise affine functions on a hierarchy of nested triangulations obtained from the initial triangulation by uniform refinement. \rv{
The triangulations are aligned with the line segments where the jumps in the coefficients take place.}

\rev{We consider a geometric multigrid V-cycle method with 6 levels to solve the discrete problems on the finest level. We generate the sequence of stiffness matrices $\A_{0:J}$, by discretizing the problems on each level of the hierarchy. The sizes of the stiffness matrices are the same for both the Poisson and the jump-1024 problems. The size of the finest-level problems is $1.64\cdot 10^6$ degrees of freedom (DoF). The size of the coarsest level problems is $1521$ DoF.
We use the standard prolongation matrices associated with the finite element spaces. The restriction matrices are transposes of the prolongation matrices.
}

The \rev{stiffness and prolongation} matrices are generated in the FE software FEniCS (version 2019.1.0) \cite{Alnaes2015,Logg2012}.
In FEniCS the stiffness matrix is assembled using all nodes of the mesh. The homogeneous Dirichlet boundary condition is then applied by setting to zero all non-diagonal elements in rows and columns which correspond to nodes on the boundary and setting to zero the corresponding elements in the right-hand side vector. We modify the stiffness matrices, the prolongation matrices and the right-hand side vector so that the Galerkin condition \eqref{eq:galerkin_condition} is satisfied. The computation is done in MATLAB 2023a. The codes for all experiments presented in this paper can be found at \url{https://github.com/vacek-petr/inVcycle}.

\rv{
Pre-smoothing and post-smoothing \rev{in the V-cycle method} are each accomplished via one iteration of the symmetric Gauss-Seidel method.
We consider the symmetric Gauss-Seidel smoother in the experiments in this paper since we are able to numerically approximate the convergence rate of the exV-cycle method in the $\A$-norm in this setting; see the discussion \Cref{subsec:ex1_relative_errERR} and \Cref{appendix}. The theoretical results stated in the paper, however, does not assume symmetry of the smoothing operators.}

\rv{We consider two variants of the coarsest-level solver: the MATLAB backslash operator and the conjugate gradient method (CG) \cite{Hestenes1952}. 
CG is stopped using a relative residual stopping criterion; i.e.,
for a chosen tolerance $\tau$ it is stopped when \mbox{$\|\vec{f}_0 - \A_0\vec{v}_{0,\inex}  \|  / \| \vec{f}_0 \| \leq \tau$}.
We consider various choices of the tolerance $\tau=2^{-i}$, $i=1,\ldots,20$.}

\rv{
We run the V-cycle methods starting with a zero initial approximation and stop when the $\A$-norm of the error is (approximately) lower than a tolerance $\theta$, i.e., $\| \vec{x} - \vec{x}^{(n)}_{\inex}\|_{\A} \leq \theta$. We consider two choices of the tolerance $\theta = 10^{-4}$ and $\theta = 10^{-11}$.
To approximate the $\A$-norm of the error on the finest level, we compute the solution using the MATLAB backslash operator.
}

\begin{figure}
\centering
\includegraphics[width=0.9\textwidth]{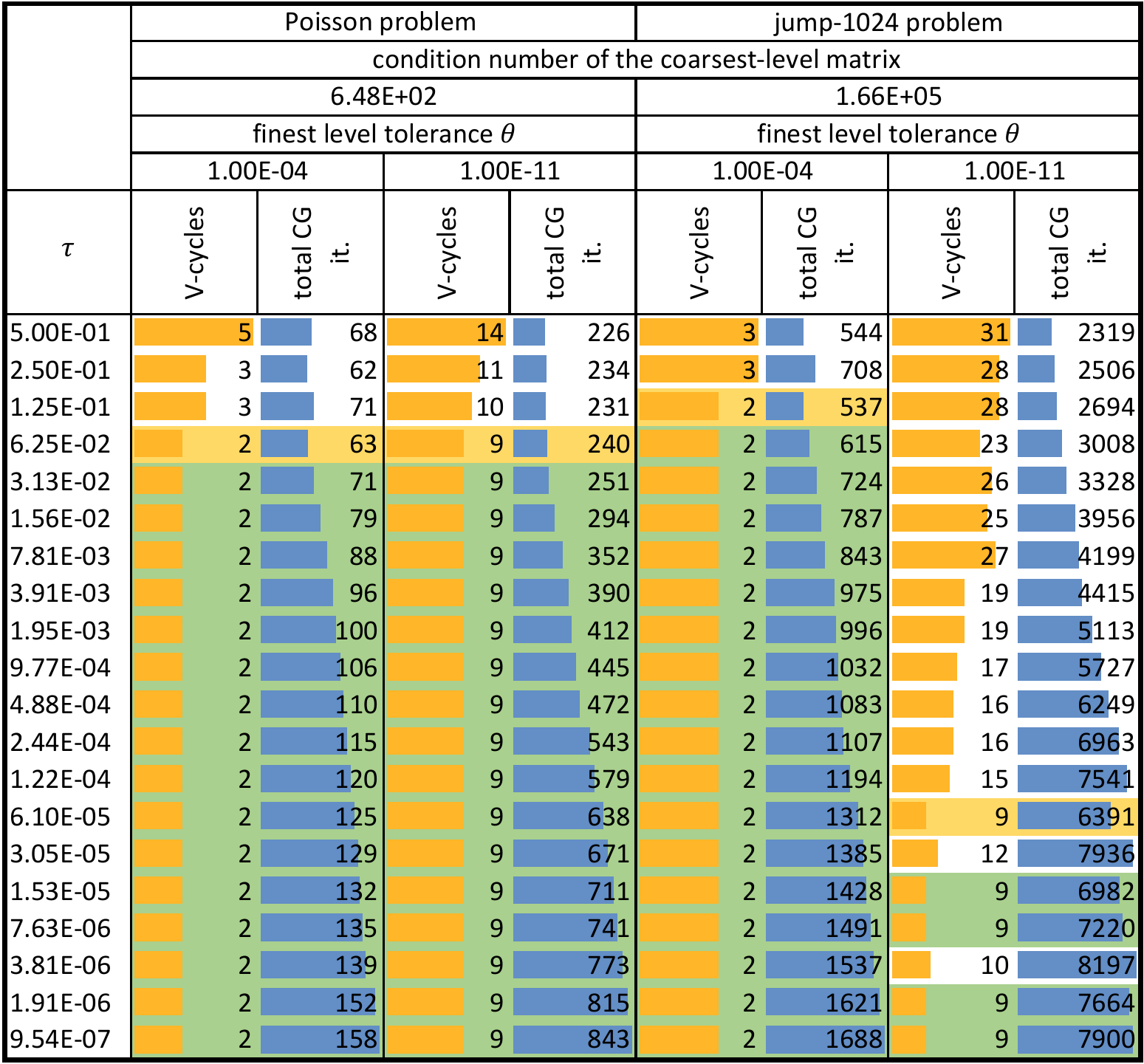}
\caption{
\rv{Comparison of inV-cycle methods with CG as the coarsest-level solver with various choices of relative residual tolerance $\tau$. The bright yellow and green color highlight variants that converge in the same number of V-cycles as the variant with MATLAB backslash operator on the coarsest-level. The bright yellow variants achieve this in the least total number of CG iterations on the coarsest-level.
}
}
\label{fig:motiv}
\end{figure}

\rv{
For both problems the variant with MATLAB backslash operator as the coarsest-level solver requires $2$ and $9$ V-cycle iterations to reach the desired finest-level accuracy $10^{-4}$ and $10^{-11}$, respectively. The results of the variants with CG as the coarsest-level solver are summarized in \Cref{fig:motiv}. 
}

\rv{
Let us first focus on the results for the Poisson problem and finest-level tolerance $\theta=10^{-4}$. The variants with CG with high coarsest-level tolerances $\tau = 2^{-i}$ ($i=1,2,3$) converge in a higher number of V-cycles than the variant with MATLAB backslash operator. The stricter the tolerance $\tau$ is the smaller the delay. The variants with tolerances $\tau = 6.25\cdot 10^{-2}$ and smaller converge in the same number of V-cycles as the method with MATLAB backslash. The variant with tolerance $\tau = 6.25\cdot 10^{-2}$ achieves this in the least total number of CG iterations on the coarsest level; this variant is in the figure highlighted by a bright yellow color. 
Using stricter tolerance than $\tau = 6.25 \cdot 10^{-2}$ is in this setting not beneficial \rev{since it does not yield a lower number of V-cycles} but it requires more computational work on the coarsest level.
We see analogous behavior for the Poisson problem and finest-level tolerance $\theta=10^{-11}$. The bright yellow highlighted variant has the same coarsest-level tolerance.
}

\rv{
Let us now focus on the results for the jump-1024 problem. The coarsest-level problem used when solving the jump-1024 problem has higher condition number than the one used for solving the Poisson problem.
The total number of coarsest-level CG iterations is for all variants significantly higher than for the corresponding variants for the Poisson problem.
We again see that the variants with high tolerances converge in a higher number of V-cycles than the variants with MATLAB backslash operator and that this delay becomes smaller for a lower coarsest-level tolerances and eventually vanishes if the tolerance is sufficiently small. It, however, does not strictly hold that \rev{lowering the tolerance results in faster converge}. This can be seen for example when comparing the variants with tolerance $\tau = 3.05\cdot10^{-5}$  and $\tau = 6.10\cdot10^{-5}$ in the setting with $\theta=10^{-11}$.
In contrast to the methods for the Poisson problem (where the values of the tolerance of the bright yellow highlighted variants are the same for the two different finest-level tolerances $\theta$) in the setting with the jump-1024 problem these values changes significantly - in order to reach the higher finest-level accuracy in the same number of V-cycles as the variant with MATLAB backslash solver the coarsest-level tolerance has to be significantly lower.
}

These experiments demonstrate that the choice of coarsest-level solver accuracy can significantly affect the convergence behavior of the V-cycle method and the overall amount of work that has to be done. This relationship is not yet well understood. 
This leads us to pose the following questions, which drive the work in this paper.
\begin{enumerate}
    \item Can we analytically describe how the accuracy of the solver on the coarsest level affects the convergence behavior of the V-cycle method?
    \item Can we define coarsest-level stopping criteria that would yield a computed V-cycle approximation ``close'' to the V-cycle approximation which would be obtained by solving the coarsest-level problems exactly?
\end{enumerate}
\section{Convergence analysis of the inV-cycle method}
\label{sec:inV-cycle_analysis}
We start by stating a few results and assumptions on the convergence of the exV-cycle method.
Let $\vec{x}^{\new}_{\ex}$ be an approximation computed by one iteration of the exV-cycle method starting with an approximation $\vec{x}^{\prev}$.
The error of the approximation $\vec{x}^{\new}_{\ex}$ can be written as the error of the previous approximation $\vec{x}^{\prev}$ times the error propagation matrix\footnote{\rv{
The error propagation matrix for a two-level exV-cycle method can be expressed as
\begin{equation*}
\matrx{E} = (\matrx{I}_1-\matrx{N}_1\matrx{A}_1)(\matrx{I}_1-\matrx{P}_1 \matrx{A}^{-1}_0 \matrx{P}^{\top}_1\matrx{A}_1)(\matrx{I}_1-\matrx{M}_1\matrx{A}_1).
\end{equation*}
A recursive expression for the error propagation matrix for an exV-cycle method with a higher number of levels can be found, e.g., in \cite[Theorem~2.4.1]{Trottenberg2001}.
}} $\matrx{E}$, i.e., 
\begin{equation*}
\vec{x} - \vec{x}^{\new}_{\ex} 
=  \matrx{E} (\vec{x}-\vec{x}^{\prev}).
\end{equation*}
We assume that the error propagation matrix $\matrx{E}$ corresponds to an operator which is a contraction with respect to the $\A$-norm, i.e., $\| \matrx{E}\|_{\A} <1$.
Proofs of this property for geometric multigrid methods can be found, e.g., in \cite{Xu1992}, \cite{Yserentant1993}.
 The contraction property implies that each iteration of the exV-cycle method reduces the $\A$-norm of the error by at least a factor $\| \matrx{E} \|_{\A}$, i.e.,
\begin{equation*}
\|\vec{x}- \vec{x}^{\new}_{\ex} \|_{\A}  \leq \| \matrx{E} \|_{\A}  \|\vec{x}-\vec{x}^{\prev}\|_{\A} \quad \forall \vec{x}^{\prev}.
\end{equation*}
We remark that this is a worst-case scenario analysis. The actual rate of convergence depends on the right-hand side and the current approximation and cannot be accurately described by a one-number characteristic.

In contrast to the exV-cycle method, the error of the approximation computed after one iteration of the inV-cycle method might not be able to be written as an error propagation matrix times the previous error. 
This is due to the fact that we consider a general solver on the coarsest level, whose application might not be able to be expressed as a matrix times vector.
To obtain insight into the convergence behavior of the inV-cycle method, we view it as a perturbation of the exV-cycle method.

Let $\vec{x}^{\new}_{\inex}$ denote the approximation computed after one iteration of the inV-cycle method starting with $\vec{x}^{\prev}$. The error of the inV-cycle approximation can be written as the error of the approximation $\vec{x}^{\new}_{\ex}$ computed after one iteration of the exV-cycle method starting with the same $\vec{x}^{\prev}$ plus the difference of the two approximations, i.e., 
\begin{equation}\label{eq:error_decompomposition}
\vec{x}-\vec{x}^{\new}_{\inex}  
= \vec{x} - \vec{x}^{\new}_{\ex}   + \vec{x}^{\new}_{\ex}- \vec{x}^{\new}_{\inex} = 
\matrx{E} (\vec{x} - \vec{x}^{\prev}) + \vec{x}^{\new}_{\ex} - \vec{x}^{\new}_{\inex}.
\end{equation}
Taking $\A$-norms on the left and right sides, using the triangle inequality and the norm of $\matrx{E}$ yields
\begin{equation}\label{eq:error_bound_with_difference}
\|  \vec{x} - \vec{x}^{\new}_{\inex} \|_{\A}
 \leq  
\| \matrx{E} \|_{\A}  \|\vec{x} - \vec{x}^{\prev} \|_{\A} + \| \vec{x}^{\new}_{\ex} - \vec{x}^{\new}_{\inex} \|_{\A} .
\end{equation}

We turn our focus to the difference $\vec{x}^{\new}_{\ex} - \vec{x}^{\new}_{\inex}$.
When applying one step of the inV-cycle method or one step of the exV-cycle method, all intermediate results $\vec{v}^{[1]}_{j}$, $j = 1, \ldots,J$, $\vec{f}_{j}$, $j = 0, \ldots,J$ are the same until the coarsest level is reached. In the exV-cycle method, the exact solution $\vec{v}_{0}$ of the problem on the coarsest level is used, while in the inV-cycle method its computed approximation $\vec{v}_{0,\inex}$ is used.
Writing down the difference $\vec{x}^{\new}_{\ex} - \vec{x}^{\new}_{\inex}$ using the individual steps in \Cref{alg:V-cycle} yields (the subscripts ``$\ex$'' and ``$\inex$'' indicate that the term corresponds to the exV-cycle method and the inV-cycle method, respectively) 
\rv{
\begin{align*}
\vec{x}^{\new}_{\ex} - \vec{x}^{\new}_{\inex} &= 
\vec{v}^{[4]}_{J,\ex} - \vec{v}^{[4]}_{J,\inex} \\
&=\vec{v}^{[3]}_{J,\ex}  +  \matrx{N}_J(\vec{f}_J  - \matrx{A}_J\vec{v}^{[3]}_{J,\ex})  -
(\vec{v}^{[3]}_{J,\inex}  +  \matrx{N}_J(\vec{f}_J  - \matrx{A}_J \vec{v}^{[3]}_{J,\inex}) )\\
& = (\matrx{I}_J - \matrx{N}_J \matrx{A}_J) 
(\vec{v}^{[3]}_{J,\ex} - \vec{v}^{[3]}_{J,\inex}) \\
& = (\matrx{I}_J - \matrx{N}_J \matrx{A}_J) 
(\vec{v}^{[1]}_{J} + \matrx{P}_J\vec{v}^{[2]}_{J-1,\ex} -(\vec{v}^{[1]}_{J} + \matrx{P}_J\vec{v}^{[2]}_{J-1,\inex})) \\
&= (\matrx{I}_J - \matrx{N}_J \matrx{A}_J) \matrx{P}_J
( \vec{v}^{[2]}_{J-1,\ex}  - \vec{v}^{[2]}_{J-1,\inex}) \\
&= (\matrx{I}_J - \matrx{N}_J \matrx{A}_J) \matrx{P}_J
( \vec{v}^{[4]}_{J-1,\ex}  - \vec{v}^{[4]}_{J-1,\inex}) \\
& = (\matrx{I}_{J} -  \matrx{N}_{J} \matrx{A}_{J})\matrx{P}_{J} \ldots (\matrx{I}_{1} -  \matrx{N}_{1} \matrx{A}_{1})\matrx{P}_{1} (\vec{v}_{0} - \vec{v}_{0,\inex}).
\end{align*}
}
Denoting by $\Sm$ the matrix
\begin{equation} \label{eq:matrixS}
    \Sm=(\matrx{I}_{J} -  \matrx{N}_{J} \matrx{A}_{J})\matrx{P}_{J} \ldots (\matrx{I}_{1} -  \matrx{N}_{1} \matrx{A}_{1})\matrx{P}_{1} \in \R^{n_J \times n_0 }
\end{equation}
gives
\begin{equation}\label{eq:difference=S_coarse_err}
\vec{x}^{\new}_{\ex} - \vec{x}^{\new}_{\inex}  = \Sm (\vec{v}_{0} - \vec{v}_{0,\inex}).
\end{equation}
We have expressed the difference of the inV-cycle and exV-cycle approximation as a matrix $\Sm$ times the error of the coarsest-level solver. The matrix $\Sm$ describes how the error is propagated to the finest level.
Let  $\|  \Sm \|_{\matrx{A}_0,\A}$ denote the norm of $\Sm$  generated by the vector norms $\| \cdot \|_{\A_0}$ and $\| \cdot \|_{\A}$, i.e.,
\begin{equation}\label{eq:norm_of_S_definition}
    \|  \Sm \|_{\matrx{A}_0,\A} =  
    \max_{\vec{v}\in\R^{n_0}, \vec{v} \neq \vec{0}}  \frac{\| \Sm\vec{v} \|_{\A}}{\| \vec{v} \|_{\matrx{A}_0}}.
\end{equation}

We derive a bound on the norm $\|  \Sm \|_{\matrx{A}_0,\A}$.
Denoting by $\matrx{S}_j$, $j=2,\ldots,J-1$, the matrix
\begin{equation*}
    \matrx{S}_j =  (\matrx{I}_{j} -  \matrx{N}_{j} \matrx{A}_{j})\matrx{P}_{j} \ldots (\matrx{I}_{1} -  \matrx{N}_{1} \matrx{A}_{1})\matrx{P}_{1} \in \R^{n_j\times n_0},
\end{equation*}
and using the definition of $\|  \matrx{S} \|_{\matrx{A}_0,\A}$ leads to
\begin{align}
\nonumber
\| \Sm \|_{\matrx{A}_0,\A} &=  \max_{\vec{v}\in\R^{n_0}, \vec{v} \neq \vec{0}} \frac{ \|  (\matrx{I}_{J} -  \matrx{N}_{J} \matrx{A}_{J})\matrx{P}_{J}\matrx{S}_{J-1} \vec{v} \|_{\A}}{\| \vec{v}\|_{\matrx{A}_0}}  \\
\nonumber
&= \max_{\vec{v}\in\R^{n_0}, \vec{v} \neq \vec{0}} 
\frac{ \|  (\matrx{I}_{J} -  \matrx{N}_{J} \matrx{A}_{J})\matrx{P}_{J}\matrx{S}_{J-1} \vec{v} \|_{\A}}
{\|  \matrx{P}_{J}\matrx{S}_{J-1} \vec{v} \|_{\A}} 
\frac{\|  \matrx{P}_{J}\matrx{S}_{J-1} \vec{v} \|_{\A}}
{\| \vec{v}\|_{\matrx{A}_0}} \\
\nonumber
& \leq \max_{\vec{v}\in\R^{n_0}, \vec{v} \neq \vec{0}} 
\| \matrx{I}_{J} -  \matrx{N}_{J} \matrx{A}_{J}\|_{\A}
\frac{\|  \matrx{P}_{J}\matrx{S}_{J-1} \vec{v} \|_{\A}}
{\| \vec{v}\|_{\matrx{A}_0}}
\\ \label{eq:galerkin_og_used}
& = 
\| \matrx{I}_{J} -  \matrx{N}_{J} \matrx{A}_{J}\|_{\A}\max_{\vec{v}\in\R^{n_0}, \vec{v} \neq \vec{0}}  
\frac{\|  \matrx{S}_{J-1} \vec{v} \|_{\matrx{A}_{J-1}}}
{\| \vec{v}\|_{\matrx{A}_0}}\\
\nonumber
& \leq \prod_{j=1}^{J} \| \matrx{I}_{j} -  \matrx{N}_{j} \matrx{A}_{j}\|_{\matrx{A}_j} 
\max_{\vec{v}\in\R^{n_0}, \vec{v} \neq \vec{0}}  
\frac{\|  \matrx{I}_{0} \vec{v} \|_{\matrx{A}_{0}}}
{\| \vec{v}\|_{\matrx{A}_0}}  \\
\nonumber
 & =  \prod_{j=1}^{J} \| \matrx{I}_{j} -  \matrx{N}_{j} \matrx{A}_{j}\|_{\matrx{A}_j},
\end{align}
where we have used the Galerkin condition \eqref{eq:galerkin_condition} to obtain \eqref{eq:galerkin_og_used}.
The monotone convergence of the post-smoothers \eqref{eq:smoothing_convergence_Anorms} in the $\A_j$-norms implies that
    $\|  \Sm \|_{\matrx{A}_0,\A}  <1$. 
If post-smoothing is not used, i.e., $\matrx{N}_j=\vec{0}$, then $\| \Sm \|_{\matrx{A}_0,\A}=1$.

The relation \eqref{eq:difference=S_coarse_err} implies
\begin{equation}\label{eq:difference_one_V_cycle_coarsest_err}
\|  \vec{x}^{\new}_{\ex} - \vec{x}^{\new}_{\inex} \|_{\A} \leq  
  \|  \matrx{S} \|_{\matrx{A}_0,\A} \| \vec{v}_{0} - \vec{v}_{0,\inex} \|_{\A_0}.
\end{equation}
Returning back to the estimate of the $\A$-norm of the error of the inV-cycle approximation, using \eqref{eq:error_bound_with_difference} and \eqref{eq:difference_one_V_cycle_coarsest_err} we have
\begin{equation}\label{eq:error_bound_with_difference_coarse}
\|  \vec{x} - \vec{x}^{\new}_{\inex} \|_{\A}
 \leq  
\| \matrx{E} \|_{\A} 
 \|\vec{x} - \vec{x}^{\prev} \|_{\A} + \|  \matrx{S} \|_{\matrx{A}_0,\A} \| \vec{v}_{0} - \vec{v}_{0,\inex} \|_{\A_0}, \quad \forall \vec{x}^{\prev}.
\end{equation}

\rv{We consider two different assumptions on the $\A_0$-norm of the error of the approximate coarsest-level solver $\| \vec{v}_{0} - \vec{v}_{0,\inex} \|_{\A_0}$:}
\rv{
\begin{itemize}
\item A \emph{relative} assumption, where the $\A_0$-norm of the error of the coarsest-level solver is less than a factor of the $\A$-norm of the error of the previous approximation on the finest level, i.e.,
there is a constant $\gamma>0$ such that
\begin{equation}\label{eq:gamma_def}
\|  \vec{v}_0- \vec{v}_{0,\inex}\|_{\matrx{A}_0}  \leq \gamma \|  \vec{x}  -  \vec{x}^{\prev} \|_{\A}, \quad \forall \vec{x}^{\prev}.
\end{equation}
\item An \emph{absolute} assumption, where the $\A_0$-norm of the error of the coarsest-level solver is less than a constant, i.e., 
there is a constant $\epsilon>0$ such that
\begin{equation}\label{eq:epsilon_def}
\|  \vec{v}_0- \vec{v}_{0,\inex}\|_{\matrx{A}_0}  \leq \epsilon , \quad \forall \vec{x}^{\prev}.
\end{equation}
\end{itemize}
}
\rv{
We first analyze the inV-cycle method under the relative assumption and then under the absolute assumption.
We comment on verification of the assumptions later in \Cref{sec:effects_rel_res,sec:new_stopping_criteria}.
}
\subsection{Relative coarsest-level accuracy}\label{subsec:analysis_relative}
Combining \eqref{eq:difference_one_V_cycle_coarsest_err} and \eqref{eq:gamma_def} yields 
an estimate on the $\A$-norm of the relative difference of the exV-cycle and inV-cycle approximations after one V-cycle iteration
\begin{equation*}
\frac{\|  \vec{x}^{\new}_{\ex} - \vec{x}^{\new}_{\inex} \|_{\A}}{\|\vec{x} - \vec{x}^{\prev} \|_{\A}} \leq  
  \|  \matrx{S} \|_{\matrx{A}_0,\A} \gamma.
\end{equation*}
 For the $\A$-norm of the error of the inV-cycle approximation, we have using \eqref{eq:error_bound_with_difference_coarse} and \eqref{eq:gamma_def}
\begin{equation}\label{eq:inV-cycle_bound_rate}
\|  \vec{x} - \vec{x}^{\new}_{\inex} \|_{\A}
 \leq  
\left(
\| \matrx{E} \|_{\A} + 
  \|  \matrx{S} \|_{\matrx{A}_0,\A} \gamma 
  \right) \|\vec{x} - \vec{x}^{\prev} \|_{\A}.
\end{equation}
Assuming that the error of the coarsest-level solver satisfies estimate \eqref{eq:gamma_def} with $\gamma$ such that 
\begin{equation*} 
\| \matrx{E} \|_{\A} +    \|  \matrx{S} \|_{\matrx{A}_0,\A} \gamma <1,
\end{equation*}
 the inV-cycle method converges and we have a bound on its convergence rate in terms of the bound on the rate of convergence of the exV-cycle method and $ \|  \matrx{S} \|_{\matrx{A}_0,\A}  \gamma $. 

We summarize the results in the following theorem.
\begin{theorem}\label{thm:oneV-cycle}
Let $\vec{x}^{\new}_{\ex}$ be the approximation of $\vec{x} = \A^{-1}\vec{b}$ computed after one iteration of the exV-cycle method with error propagation matrix $\matrx{E}$, \mbox{$\| \matrx{E}\|_{\A}<1$}, starting with an approximation $\vec{x}^{\prev}$. Let $\vec{x}^{\new}_{\inex}$ be an approximation of $\vec{x} = \A^{-1}\vec{b}$ computed after one iteration of the inV-cycle method starting with the same approximation $\vec{x}^{\prev}$, and assume the error of the coarsest-level solver $\vec{v}_0- \vec{v}_{0,\inex}$ satisfies
\begin{equation}\label{eq:thm:gamma_def}
\|  \vec{v}_0- \vec{v}_{0,\inex}\|_{\matrx{A}_0} \leq \gamma  \|  \vec{x}  -  \vec{x}^{\prev} \|_{\A},
\end{equation}
for some constant $\gamma>0$.
Then the following estimate on the $\A$-norm of the relative difference of the exV-cycle and inV-cycle approximations after one V-cycle iteration holds:
\begin{equation}\label{eq:thm:rel_approach_diff}
\frac{\| \vec{x}^{\new}_{\ex} - \vec{x}^{\new}_{\inex}\|_{\A}}{\| \vec{x} - \vec{x}^{\prev}\|_{\A}} \leq  \|  \matrx{S} \|_{\matrx{A}_0,\A}  \gamma, 
\end{equation}
where $\Sm$ is the matrix defined in \eqref{eq:matrixS}  satisfying $\| \Sm \|_{\A_0,\A}\leq1$.
Moreover, 
\begin{equation}\label{eq:thm:rel_approach_rate}
\|  \vec{x} - \vec{x}^{\new}_{\inex} \|_{\A} \leq  
\left(
\| \matrx{E} \|_{\A} + 
  \|  \matrx{S} \|_{\matrx{A}_0,\A} \gamma 
  \right) \|\vec{x} - \vec{x}^{\prev} \|_{\A},
\end{equation}
and if the error of the coarsest-level solver satisfies \eqref{eq:thm:gamma_def} with $\gamma$ such that 
\begin{equation*} 
\| \matrx{E} \|_{\A} + \|  \matrx{S} \|_{\matrx{A}_0,\A} \gamma <1,
\end{equation*}
the inV-cycle method converges.
\end{theorem}

A multigrid method is said to be uniformly convergent if there exist a bound on the rate of convergence which is independent of the number of levels and of the size of the problem on the coarsest level; see e.g., \cite{Xu1992,Yserentant1993}. If we assume that the exV-cycle method converges uniformly and the error of the coarsest-level solver in the inV-cycle method satisfies \eqref{eq:thm:gamma_def} with $\gamma$ such that $\| \matrx{E} \|_{\A} + \gamma <1$ holds and $\gamma$ is independent of the number of levels and the size of the problem on the coarsest level, inequality \eqref{eq:thm:rel_approach_rate} and the fact that $\| \Sm \|_{\A_0,\A}<1$ yield that the inV-cycle method converges uniformly.

\rv{We use the results presented in this section to discuss what may be the effect of the choice of tolerance in a relative residual coarsest-level stopping criterion  on the convergence of the V-cycle method in \Cref{sec:effects_rel_res}. We present numerical experiments testing the accuracy of the estimates \eqref{eq:thm:rel_approach_diff} and \eqref{eq:thm:rel_approach_rate} in \Cref{subsec:ex1_relative_errERR}.}

\subsection{Absolute coarsest-level accuracy}\label{subsec:analysis_absolute}
\rv{
We further focus on the analysis of the inV-cycle method under the assumption on the absolute coarsest-level accuracy \eqref{eq:epsilon_def}.
The following development is inspired by \cite[Section~4]{Vandeneshof2004}, where the authors analyze the inexact Richarson method.
}

\rv{
Let $\vec{x}^{(n)}_{\inex}$ be an approximation computed after $n$ iterations of the inV-cycle method, starting with an initial approximation $\vec{x}^{(0)}$, and assume the errors of the coarsest-level solver satisfy \eqref{eq:epsilon_def} with a constant $\epsilon>0$.
Using \eqref{eq:error_decompomposition} and \eqref{eq:difference=S_coarse_err}, the error of the $k$th approximation $\vec{x}^{(k)}_{\inex}$, $k= 1,\ldots,n$, can be written as 
\begin{equation*}
\vec{x}-\vec{x}^{(k)}_{\inex}  
= \matrx{E} (\vec{x} - \vec{x}^{(k-1)}_{\inex}) + \vec{g}^{(k)}, \quad k= 1,\ldots,n,
\end{equation*}
where  $\vec{g}^{(k)} = \matrx{S} ( \vec{v}^{(k)}_0 -  \vec{v}^{(k)}_{0,\inex} )$ and $\vec{v}^{(k)}_0 -  \vec{v}^{(k)}_{0,\inex}$ is the error of the coarsest-level solver when computing $\vec{x}^{(k)}_{\inex}$.
Let $\vec{x}^{(n)}_{\ex}$ be an approximation computed after $n$ iterations of the exV-cycle method starting with the same initial approximation $\vec{x}^{(0)}$.
The difference $\vec{x}^{(n)}_{\ex} - \vec{x}^{(n)}_{\inex}$ can be rewritten using the terms $\vec{g}^{(k)}$ as 
\begin{align*}
\vec{x}^{(n)}_{\ex} - \vec{x}^{(n)}_{\inex}  &= (\vec{x} -\vec{x}^{(n)}_{\inex}) - (\vec{x} -\vec{x}^{(n)}_{\ex})  \\
& = \matrx{E}(\vec{x} -\vec{x}^{(n-1)}_{\inex}) + \vec{g}^{(n)} - \matrx{E}^{n} (\vec{x} -\vec{x}^{(0)}) \\
& =  \matrx{E}(\matrx{E}(\vec{x} -\vec{x}^{(n-2)}_{\inex}) + \vec{g}^{(n-1)})  + \vec{g}^{(n)} - \matrx{E}^{n} (\vec{x} -\vec{x}^{(0)}) \\
& = \matrx{E}^2(\vec{x} -\vec{x}^{(n-2)}_{\inex}) + \matrx{E} \vec{g}^{(n-1)}  + \vec{g}^{(n)} - \matrx{E}^{n} (\vec{x} -\vec{x}^{(0)}) \\
&=\matrx{E}^n (\vec{x} -\vec{x}^{(0)}) + \sum^{n}_{k=1} \matrx{E}^{n-k} \vec{g}^{(k)}   -   \matrx{E}^{n} (\vec{x} -\vec{x}^{(0)})\\
& = \sum^{n}_{k=1} \matrx{E}^{n-k} \vec{g}^{(k)}.
\end{align*}
Taking the $\A$-norm of both sides, using the triangle inequality and the multiplicativity of the matrix norm $\| \cdot \|_{\A}$ we obtain  
\begin{equation}\label{eq:A-norm_difference_n_iter}
\|  \vec{x}^{(n)}_{\ex} - \vec{x}^{(n)}_{\inex}  \|_{\A}  =  \| \sum^{n}_{k=1} \matrx{E}^{n-k} \vec{g}^{(k)} \|_{\A}  
  \leq  \sum^{n}_{k=1} \| \matrx{E}\|^{n-k}_{\A}  \|\vec{g}^{(k)} \|_{\A}.
\end{equation}
Using that $\vec{g}^{(k)} = \matrx{S} ( \vec{v}^{(k)}_0 -  \vec{v}^{(k)}_{0,\inex} )$ and the norm of $\Sm$ \rev{\eqref{eq:norm_of_S_definition}} leads to
\begin{equation*} 
\|  \vec{x}^{(n)}_{\ex} - \vec{x}^{(n)}_{\inex}  \|_{\A}  
\rev{\leq
\sum^{n}_{k=1} \| \matrx{E} \|_{\A}^{n-k}    \|  \matrx{S}   (\vec{v}^{(k)}_0 -  \vec{v}^{(k)}_{0,\inex} )\|_{\A}}
\leq 
\sum^{n}_{k=1} \| \matrx{E} \|_{\A}^{n-k}    \|  \matrx{S} \|_{\matrx{A}_0,\A} \|  \vec{v}^{(k)}_0 -  \vec{v}^{(k)}_{0,\inex} \|_{\A_0}.
\end{equation*}
This bound provides information on how the accuracy of the solver on the coarsest level during the individual solves affects the $\A$-norm of the difference of the approximations $\vec{x}^{(n)}_{\ex}$ and $\vec{x}^{(n)}_{\inex}$.
}

\rv{
Using the assumption \eqref{eq:epsilon_def} and the bound for a sum of a geometric series we have
\begin{equation*}
\|  \vec{x}^{(n)}_{\ex} - \vec{x}^{(n)}_{\inex}  \|_{\A}  
\leq \sum^{n}_{k=1} \| \matrx{E} \|_{\A}^{n-k}    \|  \matrx{S} \|_{\matrx{A}_0,\A} \epsilon  < \|  \matrx{S} \|_{\matrx{A}_0,\A} \epsilon \sum^{+ \infty}_{\ell=0} \| \matrx{E} \|_{\A}^{\ell} \leq 
     \frac{\epsilon  \|  \matrx{S} \|_{\matrx{A}_0,\A}}{1-\| \matrx{E} \|_{\A}}.
\end{equation*}
Using the triangle inequality yields 
\begin{equation*}
\|  \vec{x} - \vec{x}^{(n)}_{\inex}  \|_{\A}  \leq \|  \vec{x} - \vec{x}^{(n)}_{\ex}  \|_{\A}  + \frac{\epsilon \|  \matrx{S} \|_{\matrx{A}_0,\A}}{1-\| \matrx{E} \|_{\A}};
\end{equation*}
i.e., the $\A$-norm of the error after $n$ V-cycle iterations is less than the $\A$-norm of the error of the exV-cycle approximation computed after $n$ V-cycles plus the term $\frac{\epsilon \|  \matrx{S} \|_{\matrx{A}_0,\A}}{1-\| \matrx{E} \|_{\A}}$.
}

\rv{
We summarize the results of this section in the following theorem.
\begin{theorem}\label{thm:absolute_approach}
\rv{
Let $\vec{x}^{(n)}_{\ex}$ be the approximation of $\vec{x} = \A^{-1}\vec{b}$ computed after $n$ iterations of the exV-cycle method with error propagation matrix $\matrx{E}$, \mbox{$\| \matrx{E}\|_{\A}<1$}, starting with an approximation $\vec{x}^{(0)}$. Let $\vec{x}^{(n)}_{\inex}$ be an approximation of $\vec{x} = \A^{-1}\vec{b}$ computed after $n$ iterations of the inV-cycle method, starting with the same approximation, and assume the errors of the coarsest-level solver $\vec{v}^{(k)}_0- \vec{v}^{(k)}_{0,\inex}$ satisfy
\begin{equation}\label{eq:thm:eps_def}
\|  \vec{v}^{(k)}_0- \vec{v}^{(k)}_{0,\inex}\|_{\matrx{A}_0} \leq \epsilon , \quad k = 1,\ldots,n,
\end{equation}
for a constant $\epsilon > 0$. Then the following estimate on the $A$-norm of the difference of $\vec{x}^{(n)}_{\ex}$ and $ \vec{x}^{(n)}_{\inex}$ holds:
\begin{equation}\label{eq:thm:absolute_aproach_difference_estimate}
 \|  \vec{x}^{(n)}_{\ex}  - \vec{x}^{(n)}_{\inex}  \|_{\A}\leq  \frac{\epsilon \|  \matrx{S} \|_{\matrx{A}_0,\A}}{1-\| \matrx{E} \|_{\A}},
\end{equation}
where $\Sm$ is the matrix defined in \eqref{eq:matrixS} and $\| \Sm \|_{\A_0,\A}\leq1$.
Moreover,
\begin{equation*}
\|  \vec{x} - \vec{x}^{(n)}_{\inex}  \|_{\A}  \leq \|  \vec{x} - \vec{x}^{(n)}_{\ex}  \|_{\A}  + \frac{\epsilon \|  \matrx{S} \|_{\matrx{A}_0,\A}}{1-\| \matrx{E} \|_{\A}}.
\end{equation*}
}
\end{theorem}
}

\rv{
We derive a coarsest-level stopping criteria based on these results in \Cref{sec:new_stopping_criteria} and perform numerical experiments studying the behavior of an inV-cycle method with the assumption on an absolute coarsest-level accuracy in \Cref{subsec:ex2_absolute_ERR}.}

\section{Effects of the choice of the tolerance in relative residual stopping criterion}
\label{sec:effects_rel_res}

Stopping an iterative coarsest-level solver based on the size of the relative residual is frequently done both in the literature and in practice. One chooses a tolerance $\tau$ and stops the solver when 
\begin{equation}\label{eq:relative_res_stop_crit}
    \frac{\|\vec{f}_0 - \A_0\vec{v}_{0,\inex}  \| }{\| \vec{f}_0 \|} \leq \tau.
\end{equation}
In this section we use the results from \Cref{subsec:analysis_relative} to analyze the effect of the choice of the tolerance $\tau$ on the convergence of the inV-cycle method. 
We show that if inequality \eqref{eq:relative_res_stop_crit} holds then inequality
\eqref{eq:thm:gamma_def} holds with a certain $\gamma$ depending on the tolerance $\tau$, and consequently we may use the results from \Cref{thm:oneV-cycle}. 

We start by showing that the Euclidean norm of the right-hand side on the coarsest level can be bounded by the Euclidean norm of the residual of the previous approximation on the finest level.
Rewriting $\vec{f}_0$ using the individual steps in \Cref{alg:V-cycle}, we have (note that $\vec{v}^{[0]}_j=\vec{0}$, $j=1,\ldots,J-1$)
\begin{align}\label{eq:rhs_coarse_rewrite}
\begin{split}
\vec{f}_0  & = \matrx{P}^{\top}_1 (\vec{f}_1 - \A_1 \vec{v}^{[1]}_1)  =  \matrx{P}^{\top}_1 (\vec{f}_1 - \A_1 ( \vec{v}^{[0]}_1 + \matrx{M}_1(\vec{f}_1 - \A_1 \vec{v}^{[0]}_1)) \\
& = \matrx{P}^{\top}_1 (\matrx{I}_1 - \A_1 \matrx{M}_1)\vec{f}_1  =  \prod_{j=1}^{J-1} \matrx{P}^{\top}_{j} (\matrx{I}_{j} -   \matrx{A}_{j}\matrx{M}_{j})\vec{f}_{J-1} .
\end{split}
\end{align}
The vector $\vec{f}_{J-1}$ can be expressed as
\begin{align}\label{eq:rhs_fine_minus_one_rewrite}
\begin{split}
\vec{f}_{J-1} & =  \matrx{P}^{\top}_J (\vec{b} - \A \vec{v}^{[1]}_J) 
 =  \matrx{P}^{\top}_J (\vec{b} - \A ( \vec{x}^{\prev} + \matrx{M}_J(\vec{b} - \A \vec{x}^{\prev}))) \\
 & = \matrx{P}^{\top}_J (\matrx{I}_J - \A \matrx{M}_J)(\vec{b} - \A \vec{x}^{\prev}) .
\end{split}
\end{align}
Denoting by $\matrx{T}$ the matrix
\begin{equation*}
\matrx{T} = \prod_{j=1}^{J} \matrx{P}_{j}^{\top}(\matrx{I}_{j} -   \matrx{A}_{j}\matrx{M}_{j}),
\end{equation*}
and combining \eqref{eq:rhs_coarse_rewrite} and \eqref{eq:rhs_fine_minus_one_rewrite}, we have $\vec{f}_0 = \matrx{T} \left( \vec{b} - \A \vec{x}^{\prev} \right)$. The matrix $\matrx{T}$ describes how the residual from the finest level is propagated to the coarsest level.
Based on this relation, we can estimate the Euclidean norm of $\vec{f}_0$ as
\begin{equation}
\label{eq:rhs_coarse_level_euclid_norm_bound}
\|\vec{f}_0\|
\leq   \|\matrx{T} \|  \| \vec{b} - \A \vec{x}^{\prev}  \|.
\end{equation}
The norm of $\matrx{T}$ can be bounded as
\begin{equation*}
    \| \matrx{T} \| \leq \prod^{J}_{j=1} \| \matrx{P}^{\top}_j \| \| \matrx{I}_j - \A_j \matrx{M}_j \|,
\end{equation*}
by a procedure analogous to that used in bounding the norm of $\| \Sm \|_{\A_0,\A}$; see \Cref{sec:inV-cycle_analysis}. 

\rv{
Utilizing \eqref{eq:rhs_coarse_level_euclid_norm_bound} to bound the term $\| \vec{f}_0\|$ in \eqref{eq:relative_res_stop_crit}, we obtain
\begin{equation*}
    \frac{\|\vec{f}_0 - \A_0\vec{v}_{0,\inex}  \| }{\|\matrx{T} \|  \| \vec{b} - \A \vec{x}^{\prev}  \|} \leq \tau.
\end{equation*}
Using that the Euclidean norm of the coarsest-level residual can be bounded from below by the $\A_0$-norm of the coarsest-level error as \rev{(see \Cref{apendix:inequalities})}
\begin{equation}
\label{eq:error_A_norm_by_residual_Eucl_norm_bound_upper}
\|\matrx{A}_0^{-1} \|^{-\frac{1}{2}} \|\vec{v}_0 - \vec{v}_{0,\inex}  \|_{\matrx{A}_0} 
\leq   \| \vec{f}_0 - \matrx{A}_0 \vec{v}_{0,\inex} \|,
\end{equation}
and that the Euclidean norm of the finest-level residual can be bounded from above by $\A$-norm of the error as \rev{(see \Cref{apendix:inequalities})}
\begin{equation}
\label{eq:error_A_norm_by_residual_Eucl_norm_bound_lower}
\|  \vec{b} -  \A\vec{x}^{\prev} \| 
  \leq   \| \matrx{A}\|^{\frac{1}{2}}   \| \vec{x} - \vec{x}^{\prev} \|_{\matrx{A}},
\end{equation}
we get
\begin{equation}\label{eq:rel_res_gamma}
\frac{\|\matrx{A}_0^{-1} \|^{-\frac{1}{2}} \|\vec{v}_0 - \vec{v}_{0,\inex} \|_{\A_0}}{\|\matrx{T} \|  \| \matrx{A}\|^{\frac{1}{2}}  \| \vec{x} - \vec{x}^{\prev} \|_{\matrx{A}} } \leq \tau,
\end{equation}
}
i.e., the inequality \eqref{eq:thm:gamma_def} holds with $\gamma = \tau   \| \matrx{T} \|   \| \matrx{A}\|^{\frac{1}{2}}  \| \A^{-1}_0\|^{\frac{1}{2}}$. Using the results from \Cref{thm:oneV-cycle}, we have an answer to the question of how the choice of the tolerance in the relative residual stopping criterion for the coarsest-level solver affects the convergence of the V-cycle method. 

We note that since \eqref{eq:rel_res_gamma} was derived using the estimates \eqref{eq:error_A_norm_by_residual_Eucl_norm_bound_upper}-\eqref{eq:error_A_norm_by_residual_Eucl_norm_bound_lower}, which may be a large overestimate, the resulting estimates may be loose and the actual quantities much smaller. 
We carry out numerical experiments investigating the accuracy of the estimates for the methods used in the motivating numerical experiment in \Cref{subsec:ex3_relative_relres}.

\section{Absolute coarsest-level stopping criteria}
\label{sec:new_stopping_criteria}
\rv{
In this section, we focus on the second question formulated after the motivational experiment; that is: 
\begin{quote}
``Can we define coarsest-level stopping criteria that would yield a computed V-cycle approximation ``close'' to the V-cycle approximation which would be obtained by solving the coarsest-level problems exactly?'' 
\end{quote}
}

\rv{
We present a new stopping criteria motivated by the assumption on an absolute accuracy of the coarsest-level solver and the results in \Cref{thm:absolute_approach}.
The inequality \eqref{eq:thm:eps_def} in the assumption on an absolute accuracy of the coarsest-level solver can not be directly used in practice as a coarsest-level stopping criterion since it involves the $\matrx{A}_0$-norm of the coarsest-level error, which is not available. We may, however, formulate coarsest-level stopping criteria using estimates of the $\matrx{A}_0$-norm of the error.
Let $\eta(\vec{v}^{(k)}_{0,\inex})$ be an upper bound on the $\matrx{A}_0$-norm of the error of the coarsest-level solver in the $k$th V-cycle iteration, i.e., 
\begin{equation}\label{eq:coarsest_level_eta_estimate}
\| \vec{v}^{(k)}_0 - \vec{v}^{(k)}_{0,\inex} \|_{\matrx{A}_0} \leq \eta(\vec{v}^{(k)}_{0,\inex}),\quad k=1,\ldots,n.  
\end{equation} 
We formulate a stopping criterion with a parameter $\epsilon>0$, which is chosen by the user, as
\begin{equation}\label{eq:computable_absolute_stopping_crit}
       \eta(\vec{v}^{(k)}_{0,\inex}) \leq \epsilon ,\quad k=1,\ldots , n.
\end{equation}
If \eqref{eq:computable_absolute_stopping_crit} holds then \eqref{eq:thm:eps_def} holds and from \Cref{thm:absolute_approach} we know that the $\matrx{A}$-norm of the difference of the inV-cycle and exV-cycle approximations after $n$ V-cycle iterations is bounded according to 
\begin{equation}\label{eq:rel_diff_abs_stop_E_norm}
 \|  \vec{x}^{(n)}_{\ex}  - \vec{x}^{(n)}_{\inex}  \|_{\A}\leq  \frac{\epsilon }{1-\| \matrx{E} \|_{\A}};
\end{equation}
here we have bounded $\|  \matrx{S} \|_{\matrx{A}_0,\A}$ by one from above.
We note that the accuracy of this estimate is influenced by the accuracy of the estimates \eqref{eq:coarsest_level_eta_estimate}. The term $\| \matrx{E}\|_{\A}$ is in general unknown. It is, however, included here in the form $1/(1-\| \matrx{E} \|_{\A})$. If we assume that $\| \matrx{E} \|_{\A} < \alpha$, (where, e.g., $\alpha=1/2$ or $\alpha=2/3$) we get  
\begin{equation}\label{eq:rel_diff_abs_stop}
\|  \vec{x}^{(n)}_{\ex}  - \vec{x}^{(n)}_{\inex}  \|_{\A}\leq  \frac{\epsilon}{1-\alpha}. 
\end{equation}
Due to the structure of the term $1/(1-\| \matrx{E} \|_{\A})$ this is not a significant overestimation even if the actual value of $\| \matrx{E} \|_{\A}$ is much smaller than $\alpha$. We note that assuming that $\| \matrx{E} \|_{\A} < 1/2$ or $\| \matrx{E} \|_{\A} < 2/3$ is a valid assumption for a well set up V-cycle methods.}

\rv{
The stopping criterion \eqref{eq:computable_absolute_stopping_crit} thus enable us to control the difference of the inV-cycle and exV-cycle approximations after $n$ V-cycles and consequently also the accuracy of the inV-cycle approximation.
If we want to compute an inV-cycle approximation whose $\A$-norm of the error is approximately at the level  $\theta$ (where e.g., $\theta=10^{-4}$ or $\theta=10^{-11}$) we may set $\epsilon$ as $ \epsilon =  (1-\alpha) \theta$. 
Using the triangle inequality and \eqref{eq:rel_diff_abs_stop} the $\A$-norm of the error of the inV-cycle approximation is bounded as
\begin{equation*}
\|  \vec{x} - \vec{x}^{(n)}_{\inex}  \|_{\A}  \leq \|  \vec{x} - \vec{x}^{(n)}_{\ex}  \|_{\A} 
 + \|  \vec{x}^{(n)}_{\ex} - \vec{x}^{(n)}_{\inex}  \|_{\A}  \leq \|  \vec{x} - \vec{x}^{(n)}_{\ex}  \|_{\A}  + \theta.
\end{equation*}
If we perform sufficiently many V-cycle iterations such that the $\A$-norm of the exV-cycle approximation (i.e.,  $\|  \vec{x} - \vec{x}^{(n)}_{\ex}  \|_{\A}$) would be approximately at the level of $\theta$, than the error of the inV-cycle approximation, $\|  \vec{x} - \vec{x}^{(n)}_{\inex}  \|_{\A}$,  is approximately at the level of $\theta$.}

\rv{The coarsest-level stopping criterion does not provide a finest-level stopping criterion for the inV-cycle method. We comment on a heuristic finest-level stopping indicator when discussing the results of numerical experiments in \Cref{subsec:ex2_absolute_ERR,subsec:ex4_absolute_stopping}.
}

\rv{
We further comment on the choice of the estimate $\eta$ on the $\A_0$-norm of the error on the coarsest-level.
We may use the residual based estimate on the $\matrx{A}_0$-norm of the error \eqref{eq:error_A_norm_by_residual_Eucl_norm_bound_upper}; i.e.,
\begin{equation}\label{eq:res_upper_bound}
        \eta(\vec{v}_{0,\inex}) = \| \matrx{A}^{-1}_0\|^{-\frac{1}{2}} \|  \vec{f}_0 - \matrx{A}_0 \vec{v}_{0,\inex} \|.
\end{equation}
The term $\| \matrx{A}^{-1}_0\|$, i.e., the reciprocal value of the smallest eigenvalue of $\A_0$, has to be in practical computations estimated or computed approximately.
}

\rv{
When we are using the conjugate gradient method or the preconditioned conjugate gradient method, we may use some of the upper bounds on the $\A_0$-norm of the error described e.g., in \cite{Golub2010} and the references therein, as well as in \cite{Calvetti2000,Meurant2019,Meurant2021,Meurant2023}. Most of these estimates are derived based on the interpretation of CG as a procedure for computing a Gauss quadrature approximation to a Riemann-Stieltjes integral. 
}

\rv{
We test the accuracy of estimate \eqref{eq:rel_diff_abs_stop} and the performance of the stopping criterion in numerical experiments in \Cref{subsec:ex4_absolute_stopping,}.
}

\section{Numerical experiments}
\label{sec:numerical_experiments}
\rv{
In this section we present numerical experiments illustrating some of the key results derived in this paper. We consider the same model problems and analogous V-cycle methods as in the motivating experiments in \Cref{section.motivating-experiment}. To approximate the errors on the finest and coarsest level we compute the solutions using the MATLAB backslash operator. We simulate the exV-cycle method by using MATLAB backslash operator as the solver on the coarsest level.
}

\subsection{inV-cycle method satisfying the relative coarsest-level accuracy assumption}
\label{subsec:ex1_relative_errERR}
\rv{
In this experiment, we study the behavior of the inV-cycle method with a coarsest-level solver which is stopped when the assumption on a relative coarsest-level accuracy is satisfied and examine the accuracy of the estimates presented in \Cref{thm:oneV-cycle}.
}

\rv{
We consider the same problems and analogous V-cycle methods as in the motivational experiments in \Cref{section.motivating-experiment}. The only difference is that we stop CG on the coarsest level when inequality \eqref{eq:thm:gamma_def} (approximately) holds, i.e., when 
\begin{equation*}
 \| \vec{v}_0 - \vec{v}_{0,\inex} \|_{\A_0} \leq  \gamma \| \vec{x} - \vec{x}^{\prev} \|_{\A}.
\end{equation*}
We consider three choices of the constant $\gamma$, $\gamma=0.3$, $\gamma=10^{-3}$, and $\gamma=10^{-4}$.
We run the V-cycle method starting with a zero initial approximate solution
and stop when the $\A$-norm of the error on the finest-level is (approximately) lower than $10^{-11}$.
}
\begin{figure}
\setlength\figureheight{18cm}
\setlength\figurewidth{0.99\textwidth}
\input{img/inV-cycleCGRelativeCoarsestLevelAccuracy.tikz}
\caption{
\rv{
Properties of inV-cycle methods with CG as the solver on the coarsest level, which is stopped when the assumption on the relative coarsest-level accuracy \eqref{eq:thm:gamma_def} is satisfied with $\gamma=0.3$ (\ref{line:ex1_gamma_0.3}), $\gamma=10^{-3}$  (\ref{line:ex1_gamma_1e-3}),
or $\gamma=10^{-4}$ (\ref{line:ex1_gamma_1e-4}). The dashed lines corresponds to the estimates $\| \matrx{E} \|_{\A} + \gamma$. For comparison we also include results of the exV-cycle method (\ref{line:ex1_backslash}).
}
}
\label{fig:ex1_30_11_23}
\end{figure}
\begin{figure}
\setlength\figureheight{18cm}
\setlength\figurewidth{0.99\textwidth}
\definecolor{mycolor1}{RGB}{59, 82, 139}%
\definecolor{mycolor2}{RGB}{94, 201, 98}%
\definecolor{mycolor3}{RGB}{33, 145, 140}%
\definecolor{blackk}{RGB}{68, 1, 84}%

\begin{tikzpicture}

\begin{axis}[%
width=0.45\figurewidth,
height=0.225\figureheight,
at={(0\figurewidth,0.25\figureheight)},
title = {Poisson problem, 6 levels},
xlabel={V-cycle iter.},
ylabel = $\| \vec{x} - \vec{x}^{(n)}\|_{\A} $,
xmin = 0,
xmax = 50,
yminorticks=true,
ymode = log,
ytick = {1e-13,1e-11,1e-4,1},
ymax = 1,
axis background/.style={fill=white}
]
\addplot [color=blackk, style = solid, line width=1.0pt,  mark=square, mark options={solid}, mark repeat = 3, forget plot]
  table[row sep=crcr]{%
0	0.187467821	\\
1	0.000719975	\\
2	3.33E-05	\\
3	2.55E-06	\\
4	2.40E-07	\\
5	2.60E-08	\\
6	3.10E-09	\\
7	3.89E-10	\\
8	5.03E-11	\\
9	6.66E-12	\\
10	1.01E-12	\\
11	4.90E-13	\\
12	4.76E-13	\\
13	4.76E-13	\\
14	4.76E-13	\\
15	4.75E-13	\\
16	4.74E-13	\\
17	4.73E-13	\\
18	4.73E-13	\\
19	4.74E-13	\\
20	4.72E-13	\\
21	4.73E-13	\\
22	4.75E-13	\\
23	4.75E-13	\\
24	4.72E-13	\\
25	4.74E-13	\\
26	4.74E-13	\\
27	4.76E-13	\\
28	4.75E-13	\\
29	4.75E-13	\\
30	4.76E-13	\\
31	4.75E-13	\\
32	4.72E-13	\\
33	4.72E-13	\\
34	4.75E-13	\\
35	4.75E-13	\\
36	4.75E-13	\\
37	4.76E-13	\\
38	4.76E-13	\\
39	4.75E-13	\\
40	4.73E-13	\\
41	4.73E-13	\\
42	4.73E-13	\\
43	4.73E-13	\\
44	4.76E-13	\\
45	4.74E-13	\\
46	4.73E-13	\\
47	4.72E-13	\\
48	4.74E-13	\\
49	4.78E-13	\\
50	4.73E-13	\\
};
\label{line:ex1_fix_backslash}
\addplot [color=mycolor1, style = solid, line width=1.0pt, mark=triangle, mark options={solid}, mark repeat = 3, forget plot]
  table[row sep=crcr]{%
0	0.187467821	\\
1	0.046843762	\\
2	0.013624983	\\
3	0.003576267	\\
4	0.001001646	\\
5	0.000280022	\\
6	7.62E-05	\\
7	2.15E-05	\\
8	5.90E-06	\\
9	1.68E-06	\\
10	4.68E-07	\\
11	1.32E-07	\\
12	3.35E-08	\\
13	8.18E-09	\\
14	2.39E-09	\\
15	6.42E-10	\\
16	1.87E-10	\\
17	5.20E-11	\\
18	1.34E-11	\\
19	3.59E-12	\\
20	5.70E-13	\\
21	4.51E-13	\\
22	4.47E-13	\\
23	4.45E-13	\\
24	4.43E-13	\\
25	4.42E-13	\\
26	4.42E-13	\\
27	4.41E-13	\\
28	4.40E-13	\\
29	4.40E-13	\\
30	4.40E-13	\\
31	4.39E-13	\\
32	4.39E-13	\\
33	4.39E-13	\\
34	4.38E-13	\\
35	4.38E-13	\\
36	4.38E-13	\\
37	4.38E-13	\\
38	4.38E-13	\\
39	4.37E-13	\\
40	4.37E-13	\\
41	4.37E-13	\\
42	4.37E-13	\\
43	4.37E-13	\\
44	4.37E-13	\\
45	4.37E-13	\\
46	4.37E-13	\\
47	4.37E-13	\\
48	4.37E-13	\\
49	4.37E-13	\\
50	4.37E-13	\\
};
\label{line:ex1_fix_gamma_0.3}

\addplot [color=mycolor2, style = solid, line width=1.0pt,  mark=o, mark options={solid}, mark repeat = 3, forget plot]
  table[row sep=crcr]{%
0	0.187467821	\\
1	0.000734633	\\
2	3.39E-05	\\
3	2.59E-06	\\
4	2.44E-07	\\
5	2.65E-08	\\
6	3.15E-09	\\
7	3.96E-10	\\
8	5.12E-11	\\
9	6.78E-12	\\
10	1.03E-12	\\
11	4.97E-13	\\
12	4.77E-13	\\
13	4.75E-13	\\
14	4.75E-13	\\
15	4.72E-13	\\
16	4.75E-13	\\
17	4.75E-13	\\
18	4.73E-13	\\
19	4.73E-13	\\
20	4.75E-13	\\
21	4.76E-13	\\
22	4.75E-13	\\
23	4.75E-13	\\
24	4.74E-13	\\
25	4.74E-13	\\
26	4.75E-13	\\
27	4.74E-13	\\
28	4.75E-13	\\
29	4.78E-13	\\
30	4.75E-13	\\
31	4.74E-13	\\
32	4.74E-13	\\
33	4.77E-13	\\
34	4.77E-13	\\
35	4.75E-13	\\
36	4.74E-13	\\
37	4.74E-13	\\
38	4.76E-13	\\
39	4.76E-13	\\
40	4.75E-13	\\
41	4.72E-13	\\
42	4.73E-13	\\
43	4.74E-13	\\
44	4.77E-13	\\
45	4.77E-13	\\
46	4.74E-13	\\
47	4.75E-13	\\
48	4.78E-13	\\
49	4.76E-13	\\
50	4.76E-13	\\
};
\label{line:ex1_fix_gamma_1e-3}

\addplot [color=mycolor3, style = solid, line width=1.0pt,  mark=x, mark options={solid}, mark repeat = 3, forget plot]
  table[row sep=crcr]{%
0	0.187467821	\\
1	0.000720257	\\
2	3.33E-05	\\
3	2.55E-06	\\
4	2.40E-07	\\
5	2.61E-08	\\
6	3.11E-09	\\
7	3.90E-10	\\
8	5.05E-11	\\
9	6.68E-12	\\
10	1.01E-12	\\
11	4.92E-13	\\
12	4.72E-13	\\
13	4.73E-13	\\
14	4.71E-13	\\
15	4.72E-13	\\
16	4.74E-13	\\
17	4.73E-13	\\
18	4.70E-13	\\
19	4.72E-13	\\
20	4.73E-13	\\
21	4.74E-13	\\
22	4.74E-13	\\
23	4.73E-13	\\
24	4.74E-13	\\
25	4.75E-13	\\
26	4.72E-13	\\
27	4.73E-13	\\
28	4.71E-13	\\
29	4.72E-13	\\
30	4.72E-13	\\
31	4.71E-13	\\
32	4.71E-13	\\
33	4.71E-13	\\
34	4.73E-13	\\
35	4.74E-13	\\
36	4.74E-13	\\
37	4.74E-13	\\
38	4.75E-13	\\
39	4.77E-13	\\
40	4.76E-13	\\
41	4.75E-13	\\
42	4.73E-13	\\
43	4.72E-13	\\
44	4.73E-13	\\
45	4.75E-13	\\
46	4.72E-13	\\
47	4.73E-13	\\
48	4.73E-13	\\
49	4.74E-13	\\
50	4.76E-13	\\
};
\label{line:ex1_fix_gamma_1e-4}
\end{axis}

\begin{axis}[%
width=0.45\figurewidth,
height=0.225\figureheight,
at={(0.5\figurewidth,0.25\figureheight)},
title = {jump-1024 problem, 6 levels},
xlabel={V-cycle iter.},
ylabel = $\| \vec{x} - \vec{x}^{(n)}\|_{\A} $,
xmin = 0,
xmax = 50,
yminorticks=true,
ymode = log,
ytick = {1e-13,1e-11,1e-4,1},
axis background/.style={fill=white}
]
\addplot [color=blackk, style = solid, line width=1.0pt,  mark=square, mark options={solid}, mark repeat = 3, forget plot]
  table[row sep=crcr]{%
0	0.066698707	\\
1	0.000700793	\\
2	3.51E-05	\\
3	2.84E-06	\\
4	2.80E-07	\\
5	3.12E-08	\\
6	3.76E-09	\\
7	4.75E-10	\\
8	6.16E-11	\\
9	8.14E-12	\\
10	1.09E-12	\\
11	1.60E-13	\\
12	6.51E-14	\\
13	6.19E-14	\\
14	6.16E-14	\\
15	6.15E-14	\\
16	6.17E-14	\\
17	6.16E-14	\\
18	6.18E-14	\\
19	6.14E-14	\\
20	6.17E-14	\\
21	6.19E-14	\\
22	6.17E-14	\\
23	6.16E-14	\\
24	6.17E-14	\\
25	6.16E-14	\\
26	6.16E-14	\\
27	6.16E-14	\\
28	6.16E-14	\\
29	6.15E-14	\\
30	6.13E-14	\\
31	6.17E-14	\\
32	6.13E-14	\\
33	6.13E-14	\\
34	6.14E-14	\\
35	6.11E-14	\\
36	6.12E-14	\\
37	6.18E-14	\\
38	6.15E-14	\\
39	6.17E-14	\\
40	6.15E-14	\\
41	6.11E-14	\\
42	6.14E-14	\\
43	6.15E-14	\\
44	6.15E-14	\\
45	6.16E-14	\\
46	6.13E-14	\\
47	6.10E-14	\\
48	6.19E-14	\\
49	6.14E-14	\\
50	6.18E-14	\\
};
\addplot [color=mycolor1, style = solid, line width=1.0pt, mark=triangle, mark options={solid}, mark repeat = 3, forget plot]
  table[row sep=crcr]{%
0	0.066698707	\\
1	0.016881125	\\
2	0.004731293	\\
3	0.001246246	\\
4	0.000329572	\\
5	9.05E-05	\\
6	2.61E-05	\\
7	8.81E-06	\\
8	4.23E-06	\\
9	2.36E-06	\\
10	1.46E-06	\\
11	8.90E-07	\\
12	5.55E-07	\\
13	3.42E-07	\\
14	2.13E-07	\\
15	1.32E-07	\\
16	8.14E-08	\\
17	5.04E-08	\\
18	3.12E-08	\\
19	1.93E-08	\\
20	1.20E-08	\\
21	7.40E-09	\\
22	4.58E-09	\\
23	2.83E-09	\\
24	1.75E-09	\\
25	1.08E-09	\\
26	6.68E-10	\\
27	4.11E-10	\\
28	2.53E-10	\\
29	1.56E-10	\\
30	9.66E-11	\\
31	5.96E-11	\\
32	3.65E-11	\\
33	2.27E-11	\\
34	1.40E-11	\\
35	8.68E-12	\\
36	5.35E-12	\\
37	3.30E-12	\\
38	2.03E-12	\\
39	1.24E-12	\\
40	7.56E-13	\\
41	4.59E-13	\\
42	2.74E-13	\\
43	1.67E-13	\\
44	1.04E-13	\\
45	6.98E-14	\\
46	5.29E-14	\\
47	4.77E-14	\\
48	4.64E-14	\\
49	4.55E-14	\\
50	4.54E-14	\\
};
\addplot [color=mycolor2, style = solid, line width=1.0pt,  mark=o, mark options={solid}, mark repeat = 3, forget plot]
  table[row sep=crcr]{%
0	0.066698707	\\
1	0.000702526	\\
2	3.54E-05	\\
3	3.00E-06	\\
4	6.08E-07	\\
5	3.33E-07	\\
6	2.05E-07	\\
7	1.27E-07	\\
8	7.81E-08	\\
9	4.83E-08	\\
10	2.98E-08	\\
11	1.84E-08	\\
12	1.14E-08	\\
13	7.03E-09	\\
14	4.34E-09	\\
15	2.68E-09	\\
16	1.66E-09	\\
17	1.02E-09	\\
18	6.32E-10	\\
19	3.90E-10	\\
20	2.41E-10	\\
21	1.49E-10	\\
22	9.20E-11	\\
23	5.68E-11	\\
24	3.51E-11	\\
25	2.17E-11	\\
26	1.34E-11	\\
27	8.27E-12	\\
28	5.11E-12	\\
29	3.15E-12	\\
30	1.95E-12	\\
31	1.20E-12	\\
32	7.42E-13	\\
33	4.59E-13	\\
34	2.86E-13	\\
35	1.81E-13	\\
36	1.20E-13	\\
37	8.74E-14	\\
38	7.16E-14	\\
39	6.52E-14	\\
40	6.29E-14	\\
41	6.17E-14	\\
42	6.16E-14	\\
43	6.16E-14	\\
44	6.13E-14	\\
45	6.13E-14	\\
46	6.14E-14	\\
47	6.16E-14	\\
48	6.16E-14	\\
49	6.15E-14	\\
50	6.21E-14	\\
};

\addplot [color=mycolor3, style = solid, line width=1.0pt,  mark=x, mark options={solid}, mark repeat = 3, forget plot]
  table[row sep=crcr]{%
0	0.066698707	\\
1	0.000701153	\\
2	3.51E-05	\\
3	2.84E-06	\\
4	2.80E-07	\\
5	3.12E-08	\\
6	3.76E-09	\\
7	4.88E-10	\\
8	9.43E-11	\\
9	4.49E-11	\\
10	2.73E-11	\\
11	1.69E-11	\\
12	1.04E-11	\\
13	6.43E-12	\\
14	3.98E-12	\\
15	2.46E-12	\\
16	1.52E-12	\\
17	9.42E-13	\\
18	5.86E-13	\\
19	3.67E-13	\\
20	2.33E-13	\\
21	1.54E-13	\\
22	1.08E-13	\\
23	8.32E-14	\\
24	7.10E-14	\\
25	6.57E-14	\\
26	6.35E-14	\\
27	6.27E-14	\\
28	6.22E-14	\\
29	6.20E-14	\\
30	6.15E-14	\\
31	6.15E-14	\\
32	6.17E-14	\\
33	6.16E-14	\\
34	6.13E-14	\\
35	6.13E-14	\\
36	6.15E-14	\\
37	6.11E-14	\\
38	6.16E-14	\\
39	6.18E-14	\\
40	6.18E-14	\\
41	6.19E-14	\\
42	6.17E-14	\\
43	6.17E-14	\\
44	6.14E-14	\\
45	6.13E-14	\\
46	6.12E-14	\\
47	6.16E-14	\\
48	6.14E-14	\\
49	6.12E-14	\\
50	6.13E-14	\\
};
\end{axis}

\end{tikzpicture}%
\caption{
\rv{
$\A$-norm of the error of the inV-cycle methods with CG as the solver on the coarsest level, which is stopped when the assumption on the relative coarsest-level accuracy \eqref{eq:thm:gamma_def} is satisfied with $\gamma=0.3$ (\ref{line:ex1_fix_gamma_0.3}), $\gamma=10^{-3}$  (\ref{line:ex1_fix_gamma_1e-3}),
or $\gamma=10^{-4}$ (\ref{line:ex1_fix_gamma_1e-4}). For comparison we also include the $\A$-norm of the error of the exV-cycle method (\ref{line:ex1_fix_backslash}). \rev{Every third point is marked.}
}
}
\label{fig:ex1_fix_19_12_23}
\end{figure}

\rv{
The results are summarized in \Cref{fig:ex1_30_11_23}.
After each V-cycle iteration we compute the $\A$-norms of the relative difference of the exV-cycle and inV-cycle approximations after one V-cycle iteration, i.e., 
\begin{equation}\label{eq:ex1_reldiff}
\frac{\| \vec{x}^{\new}_{\ex} - \vec{x}^{\new}_{\inex}\|_{\A} }{\| \vec{x} - \vec{x}^{\prev}\|_{\A}},
\end{equation}
for $\vec{x}^{\prev}=\vec{x}^{(k)}_{\inex}$, $k=0,1,\ldots$ . According to the estimate \eqref{eq:thm:rel_approach_diff} from \Cref{thm:oneV-cycle}, the relative difference \eqref{eq:ex1_reldiff} should be less than $ \gamma \| \matrx{S} \|_{\matrx{A}_0,\A}$, where $\| \matrx{S} \|_{\matrx{A}_0,\A}\leq 1$. Looking at the results we see that all values \eqref{eq:ex1_reldiff} are slightly less than $\gamma$ besides the ones computed after the last few V-cycle iterations of the variants with $\gamma = 10^{-4}$.
We strongly believe that these outlier are caused by the effects of finite precision arithmetic.
Dividing the computed values \eqref{eq:ex1_reldiff} (besides the mentioned outliers) by $\gamma$ and finding the maximum we get a lower bound on $ \| \matrx{S} \|_{\matrx{A}_0,\A}$, which is $0.95$ and $0.97$ for the variant with the Poisson and the jump-1024 problem, respectively.
}

\rv{
We also compute the convergence rate in the $\A$-norm, after each V-cycle iteration, i.e., 
\begin{equation} \label{eq:ex1_errAnorm_est}
\frac{\|  \vec{x} - \vec{x}^{(n)} \|_{\A} }{\|\vec{x} - \vec{x}^{(n-1)} \|_{\A}}, \quad n = 1,2,\ldots, \quad.
\end{equation}
According to the estimate \eqref{eq:thm:rel_approach_rate}, the convergence rate \eqref{eq:ex1_errAnorm_est} is bounded by $ \| \matrx{E} \|_{\A} + \gamma$; we have used that $\| \matrx{S} \|_{\matrx{A}_0,\A}\leq 1$.
We approximate the term $\| \matrx{E} \|_{\A}$  by a procedure described in \Cref{appendix}. It is approximately $0.15$ and $0.62$ for the variant with the Poisson and the jump-1024 problem, respectively.
Looking at the results we see that all the computed values of \eqref{eq:ex1_errAnorm_est} are less than the corresponding bounds.
}

\rv{
Let us first comment on the result for the Poisson problem.
The convergence rates of the variants with $\gamma = 10^{-3}$ and $\gamma = 10^{-4}$ are approximately the same as the convergence rate of the exV-cycle method.
The rates are significantly lower than its bounds in the first few V-cycle iterations, but they gradually deteriorate to approximately the value of the bound in the last V-cycle iterations.
The convergence rate of the variant with $\gamma = 0.3$ is approximately constant $0.3$. Here we don't see the usual deterioration of the convergence rate after the first V-cycle iterations. The bound for this variant is approximately $0.45$.
}

\rv{
Let us focus on the results for the jump-1024 problem.
The convergence rate of the exV-cycle method doesn't deteriorate to the value of its approximate bound $0.62$, but it stays under $0.15$. This is an interesting behaviour since $0.15$ is approximately the value of the bound on the rate of convergence of the exV-cycle method for the Poisson problem.
The convergence rates of the variants with $\gamma = 10^{-3}$, $\gamma = 10^{-4}$, are in the first V-cycle iterations approximately the same as the rate of the exV-cycle method. They, however, eventually deteriorate to the expected bounds. The deterioration happens sooner for the variant with $\gamma = 10^{-3}$.
}

\rv{
The convergence rate of the variant with $\gamma = 0.3$ is approximately $0.3$ in the first few iterations then it deteriorates to $0.62$. This is another interesting behaviour since $0.62$ is the value of the bound on the convergence rate of the exV-cycle method. The bound on the convergence rate of the inV-cycle method with $\gamma = 0.3$ is $0.92$.
}

\rv{
In these experiments we see that the estimate of the rate of convergence of the inV-cycle method with the assumption on a relative coarsest level accuracy is an accurate estimate of the worst-case convergence rate if $\gamma$ is smaller than $\| \matrx{E} \|_{\A}$.
}

\rv{
We also plot the $\A$-norm of the error and the number of CG iterations on the coarsest level. We see that the number of CG iterations performed in the variants with the jump-1024 problem is significantly higher than in the variants with the Poisson problem.
}

\rv{
To find out whether the inV-cycle methods reach the same level of attainable accuracy as the exV-cycle methods, we perform an experiment, where we stop the V-cycle method on the finest level after $50$ V-cycle iterations. The results are summarized in \Cref{fig:ex1_fix_19_12_23}. We see that the considered inV-cycle methods reach the same level of attainable accuracy as the exV-cycle methods.
}

\subsection{Accuracy of the estimates for inV-cycle methods with a relative residual coarsest-level stopping criterion}
\label{subsec:ex3_relative_relres}
\rv{
In this experiment, we study the accuracy of the results for a inV-cycle methods with a relative residual coarsest-level stopping criterion discussed in \Cref{sec:effects_rel_res}.
}

\rv{
We consider the same problems and analogous V-cycle methods as in the motivational experiments in \Cref{section.motivating-experiment}. We stop CG on the coarsest level using the relative residual stopping criterion \eqref{eq:relative_res_stop_crit}, i.e., when 
\begin{equation*}
\frac{\| \vec{f}_0 - \A_0 \vec{v}_{0,\inex} \| }{\| \vec{f}_0 \|}\leq  \tau,
\end{equation*}
and choose $\tau = 10^{-4}   \| \matrx{T} \|^{-1}   \| \matrx{A}\|^{-\frac{1}{2}}  \| \matrx{A}^{-1}_0\|^{-\frac{1}{2}}$.
We approximate the terms $\| \matrx{T} \|$, $\| \matrx{A}\|$,  $\| \matrx{A}^{-1}_0\|$ using MATLAB function \texttt{eigs}.
}

\rv{
We run the V-cycle method starting with a zero initial approximate solution
and stop when the $\A$-norm of the error on the finest level is (approximately) lower than $10^{-11}$.
In order to find out whether the results are substantially  affected by the use of the finite precision arithmetic, we run the computation both in the standard MATLAB double precision and also in a simulated quad precision using the Advanpix toolbox \cite{advanpix}.
}

\rv{
After each V-cycle iteration we compute the $\A$-norm of the relative difference of the exV-cycle and inV-cycle approximations after one V-cycle iteration \eqref{eq:ex1_reldiff}.
The V-cycle methods for both problems reach the desired accuracy in $9$ V-cycle iterations. The results are summarized in \Cref{fig:comparison_stop_crit}.
According to the discussion in \Cref{sec:effects_rel_res} 
the relative difference \eqref{eq:ex1_reldiff} should be less than 
\begin{equation*}
 \tau   \| \matrx{T} \|   \| \matrx{A}\|^{\frac{1}{2}}  \| \A^{-1}_0\|^{\frac{1}{2}} \|  \matrx{S}\|_{\A_0,\A}    .
\end{equation*}
Bounding $\|  \matrx{S}\|_{\A_0,\A}$ by one from above and considering our choice of $\tau$, we get that the relative difference \eqref{eq:ex1_reldiff} should be less than $10^{-4}$. We see that this is true for all of the computed values. The computed values are however significantly smaller than the estimate. This may be a consequence of the usage of the estimates \eqref{eq:error_A_norm_by_residual_Eucl_norm_bound_upper} and \eqref{eq:error_A_norm_by_residual_Eucl_norm_bound_lower} in the derivation of the estimates in \Cref{sec:effects_rel_res}.
}

\rv{
We see that the relative difference \eqref{eq:ex1_reldiff} for the variant computed in double precision starts increasing after the $5$th  V-cycle iterations, whereas the relative difference for the variant computed in the simulated quad precision stay approximately at the same level.
We thus strongly believe that the increase of the values computed in double is caused by the use of the finite precision arithmetic. 
}
\begin{figure}
\setlength\figureheight{13.5cm}
\setlength\figurewidth{0.99\textwidth}
\definecolor{mycolor1}{RGB}{33, 145, 140}%
\definecolor{mycolor2}{RGB}{94, 201, 98}%

\begin{tikzpicture}
\begin{axis}[%
width=0.45\figurewidth,
height=0.3\figureheight,
at={(0\figurewidth,0.66\figureheight)},
title = {Poisson problem, 6 levels},
xlabel={V-cycle iter.},
ylabel = $\frac{\| \vec{x}^{\new}_{\ex} - \vec{x}^{\new}_{\inex}\|_{\A} }{\| \vec{x} - \vec{x}^{\prev}\|_{\A}}$,
xmin = 0,
xmax = 9,
ymode = log,
ytick = {1e-4,1e-6,1e-8},
ymin = 1e-10,
ymax = 1e-3,
yminorticks=true,
axis background/.style={fill=white}
]
\addplot [color=mycolor2, style = solid, line width=1.0pt, mark=o, mark options={solid}, forget plot]
  table[row sep=crcr]{%
1	2.78E-08	\\
2	3.74E-08	\\
3	3.60E-08	\\
4	4.36E-08	\\
5	5.18E-08	\\
6	3.18E-07	\\
7	2.52E-06	\\
8	1.21E-05	\\
9	3.78E-05	\\
};
\label{line:ex3_relres}
\addplot [color=mycolor1, style = solid, line width=1.0pt, mark=x, mark options={solid}, forget plot]
  table[row sep=crcr]{%
1	2.78E-08	\\
2	3.74E-08	\\
3	3.60E-08	\\
4	4.35E-08	\\
5	3.79E-08	\\
6	4.04E-08	\\
7	3.51E-08	\\
8	3.07E-08	\\
9	2.73E-08	\\
};
\label{line:ex3_relres_quad}
\end{axis}

\begin{axis}[%
width=0.45\figurewidth,
height=0.3\figureheight,
at={(0.5\figurewidth,0.66\figureheight)},
title = {jump-1024 problem, 6 levels},
xlabel={V-cycle iter.},
ylabel = $\frac{\| \vec{x}^{\new}_{\ex} - \vec{x}^{\new}_{\inex}\|_{\A} }{\| \vec{x} - \vec{x}^{\prev}\|_{\A}}$,
xmin = 0,
xmax = 9,
ytick = {1e-4,1e-6,1e-8},
ymin = 1e-10,
ymax = 1e-3,
ymode = log,
yminorticks=true,
axis background/.style={fill=white}
]
\addplot [color=mycolor2, style = solid, line width=1.0pt, mark=o, mark options={solid}, forget plot]
  table[row sep=crcr]{%
1	4.23E-10	\\
2	2.26E-09	\\
3	1.37E-09	\\
4	1.60E-09	\\
5	5.69E-09	\\
6	4.84E-08	\\
7	3.44E-07	\\
8	1.01E-06	\\
9	2.87E-06	\\
};
\addplot [color=mycolor1, style = solid, line width=1.0pt, mark=x, mark options={solid}, forget plot]
  table[row sep=crcr]{%
1	5.40E-10	\\
2	2.46E-09	\\
3	1.24E-09	\\
4	1.62E-09	\\
5	2.06E-09	\\
6	1.77E-09	\\
7	1.48E-09	\\
8	1.23E-09	\\
9	1.09E-09	\\
};
\end{axis}

\end{tikzpicture}%
\caption{
Testing accuracy of the estimate discussed in \Cref{sec:effects_rel_res}. We consider the V-cycle method with CG as the solver on the coarsest level. CG is stopped using the relative residual stopping criterion \eqref{eq:relative_res_stop_crit} with $\tau = 10^{-4}   \| \matrx{T} \|^{-1}   \| \matrx{A}\|^{-\frac{1}{2}}  \| \matrx{A}^{-1}_0\|^{-\frac{1}{2}}$. The computation is done in standard MATLAB double precision (\ref{line:ex3_relres}) and in simulated quad precision using the Advapix toolbox (\ref{line:ex3_relres_quad}).
}
\label{fig:comparison_stop_crit}
\end{figure}

\subsection{inV-cycle method satisfying the absolute coarsest-level accuracy assumption}
\label{subsec:ex2_absolute_ERR}
\rv{
In this experiment we study the behavior of the inV-cycle method with a coarsest-level solver that is stopped when the assumption on an absolute coarsest-level accuracy is satisfied and examine the accuracy of estimates presented in \Cref{thm:absolute_approach}.
}

\rv{
We consider the same problems and analogous V-cycle methods as in the motivational experiments in \Cref{section.motivating-experiment}. The only difference is that we stop CG on the coarsest level when inequality \eqref{eq:thm:eps_def} (approximately) holds, i.e., when 
\begin{equation*}
\| \vec{v}_0 - \vec{v}_{0,\inex} \|_{\A_0} \leq \epsilon . 
\end{equation*}
We choose $\epsilon =\theta( 1 - \| \matrx{E} \|_{\A} )$, where $\theta=10^{-4}$ or $\theta=10^{-11}$. We approximate $\| \matrx{E} \|_{\A}$ as in the experiments in \Cref{subsec:ex1_relative_errERR}.
We run the V-cycle method starting with a zero initial approximate solution
and stop after $15$ V-cycle iterations.
}

\begin{figure}
\setlength\figureheight{13.5cm}
\setlength\figurewidth{0.99\textwidth}
\definecolor{mycolor1}{RGB}{33, 145, 140}%
\definecolor{mycolor2}{RGB}{33, 145, 140}%
\definecolor{mycolor8}{RGB}{68, 1, 84}%

\begin{tikzpicture}

\begin{axis}[%
width=0.45\figurewidth,
height=0.3\figureheight,
at={(0\figurewidth,0.66\figureheight)},
title = {Poisson problem, 6 levels},
xlabel={V-cycle iter.},
ylabel = $\| \vec{x}^{(n)}_{\ex} - \vec{x}^{(n)}_{\inex}\|_{\A}$,
xmin = 0,
xmax = 15,
xtick = {0,2,10,15,20},
ytick = {1e-4,1e-11},
ymax = 1e-4,
yminorticks=true,
ymode = log,
axis background/.style={fill=white}
]

\addplot[
color=mycolor1,
style = solid,
line width=1.0pt,
mark options={solid},
forget plot
]
table[row sep=crcr]{%
1	2.97E-12	\\
2	2.47E-12	\\
3	2.51E-12	\\
4	2.40E-12	\\
5	2.79E-12	\\
6	2.95E-12	\\
7	3.29E-12	\\
8	3.56E-12	\\
9	3.60E-12	\\
10	3.58E-12	\\
11	3.39E-12	\\
12	3.21E-12	\\
13	3.07E-12	\\
14	2.95E-12	\\
15	2.84E-12	\\
16	2.74E-12	\\
17	2.64E-12	\\
18	2.55E-12	\\
19	2.46E-12	\\
20	2.38E-12	\\
};
\label{line:ex2_delta_1e-11} 
\addplot[
color=mycolor2,
style = densely dashed,
line width=1.0pt,
mark options={solid},
forget plot
]
table[row sep=crcr]{%
1	3.20E-05	\\
2	2.97E-05	\\
3	2.79E-05	\\
4	2.29E-05	\\
5	1.95E-05	\\
6	1.70E-05	\\
7	1.52E-05	\\
8	1.38E-05	\\
9	1.27E-05	\\
10	1.17E-05	\\
11	1.10E-05	\\
12	1.03E-05	\\
13	9.73E-06	\\
14	9.22E-06	\\
15	8.76E-06	\\
16	8.35E-06	\\
17	7.97E-06	\\
18	7.63E-06	\\
19	7.31E-06	\\
20	7.01E-06	\\
};
\label{line:ex2_delta_1e-4}
\end{axis}

\begin{axis}[%
width=0.45\figurewidth,
height=0.3\figureheight,
at={(0.5\figurewidth,0.66\figureheight)},
title = {jump-1024 problem, 6 levels},
xlabel={V-cycle iter.},
ylabel = $\| \vec{x}^{(n)}_{\ex} - \vec{x}^{(n)}_{\inex}\|_{\A}$,
xmin = 0,
xmax = 15,
ymax = 1e-4,
xtick = {0,2,9,15},
ytick = {1e-4,1e-11},
yminorticks=true,
ymode = log,
axis background/.style={fill=white}
]
\addplot[
color=mycolor1,
style = solid,
line width=1.0pt,
mark options={solid},
forget plot
]
table[row sep=crcr]{%
1	2.86E-12	\\
2	6.48E-12	\\
3	3.39E-12	\\
4	3.21E-12	\\
5	2.71E-12	\\
6	5.15E-12	\\
7	5.21E-12	\\
8	6.11E-12	\\
9	7.16E-12	\\
10	7.21E-12	\\
11	6.15E-12	\\
12	5.27E-12	\\
13	4.57E-12	\\
14	4.00E-12	\\
15	3.51E-12	\\
16	3.09E-12	\\
17	2.72E-12	\\
18	2.40E-12	\\
19	2.12E-12	\\
20	1.88E-12	\\
}; 
\addplot[
color=mycolor2,
style = densely dashed,
line width=1.0pt,
mark options={solid},
forget plot
]
table[row sep=crcr]{%
1	3.55E-05	\\
2	4.27E-05	\\
3	5.33E-05	\\
4	4.28E-05	\\
5	3.57E-05	\\
6	3.07E-05	\\
7	2.68E-05	\\
8	2.34E-05	\\
9	2.06E-05	\\
10	1.81E-05	\\
11	1.59E-05	\\
12	1.39E-05	\\
13	1.23E-05	\\
14	1.08E-05	\\
15	9.47E-06	\\
16	8.33E-06	\\
17	7.33E-06	\\
18	6.44E-06	\\
19	5.67E-06	\\
20	4.99E-06	\\
};
\end{axis}

\begin{axis}[%
width=0.45\figurewidth,
height=0.3\figureheight,
at={(0\figurewidth,0.33\figureheight)},
title = {Poisson problem, 6 levels},
xlabel={V-cycle iter.},
ylabel = $\| \vec{x} - \vec{x}^{(n)}\|_{\A}$,
xmin = 0,
xmax = 15,
ytick = {1,1.0E-4,1.0E-11,1.0E-13},
ymax = 1,
xtick = {0,2,10,15,20,30},
yminorticks = true,
ymode = log,
axis background/.style = {fill=white}
]

\addplot[
color=mycolor1,
style = solid,
line width=1.0pt,
mark options={solid},
forget plot
]
table[row sep=crcr]{%
0	0.187467821	\\
1	0.000719975	\\
2	3.33E-05	\\
3	2.55E-06	\\
4	2.40E-07	\\
5	2.60E-08	\\
6	3.10E-09	\\
7	3.89E-10	\\
8	5.10E-11	\\
9	1.02E-11	\\
10	7.75E-12	\\
11	7.53E-12	\\
12	7.41E-12	\\
13	7.32E-12	\\
14	7.24E-12	\\
15	7.16E-12	\\
16	7.08E-12	\\
17	7.01E-12	\\
18	6.94E-12	\\
19	6.87E-12	\\
20	6.81E-12	\\
};
\addplot[
color=mycolor2,
style = densely dashed,
line width=1.0pt,
mark options={solid},
forget plot
]
table[row sep=crcr]{%
0	0.187467821	\\
1	0.000723124	\\
2	8.18E-05	\\
3	7.05E-05	\\
4	6.60E-05	\\
5	6.28E-05	\\
6	6.02E-05	\\
7	5.80E-05	\\
8	5.61E-05	\\
9	5.44E-05	\\
10	5.29E-05	\\
11	5.15E-05	\\
12	5.03E-05	\\
13	4.91E-05	\\
14	4.80E-05	\\
15	4.70E-05	\\
16	4.61E-05	\\
17	4.52E-05	\\
18	4.44E-05	\\
19	4.37E-05	\\
20	4.29E-05	\\
};
\addplot[
color=mycolor8,
style = solid,
line width=1.0pt,
mark options={solid},
forget plot
]
table[row sep=crcr]{%
0	0.187467821	\\
1	0.000719975	\\
2	3.33E-05	\\
3	2.55E-06	\\
4	2.40E-07	\\
5	2.60E-08	\\
6	3.10E-09	\\
7	3.89E-10	\\
8	5.03E-11	\\
9	6.66E-12	\\
10	1.01E-12	\\
11	4.90E-13	\\
12	4.76E-13	\\
13	4.76E-13	\\
14	4.76E-13	\\
15	4.75E-13	\\
16	4.74E-13	\\
17	4.73E-13	\\
18	4.73E-13	\\
19	4.74E-13	\\
20	4.72E-13	\\
};
\label{line:ex2_backslash}
\end{axis}

\begin{axis}[%
width=0.45\figurewidth,
height=0.3\figureheight,
at={(0\figurewidth,0\figureheight)},
title = {Poisson problem, 6 levels},
xlabel={V-cycle iter.},
ylabel = number of CG iter.,
xtick = {0,2,10,15,20,30},
xmin = 0,
xmax = 15,
yminorticks=true,
axis background/.style={fill=white}
]
\addplot[
color=mycolor1,
style = solid,
line width=1.0pt,
mark options={solid},
forget plot
]
table[row sep=crcr]{%
1	110	\\
2	111	\\
3	103	\\
4	89	\\
5	73	\\
6	56	\\
7	27	\\
8	10	\\
9	3	\\
10	0	\\
11	0	\\
12	0	\\
13	0	\\
14	0	\\
15	0	\\
16	0	\\
17	0	\\
18	0	\\
19	0	\\
20	0	\\
};
\addplot[
color=mycolor2,
style = densely dashed,
line width=1.0pt,
mark options={solid},
forget plot
]
table[row sep=crcr]{%
1	42	\\
2	27	\\
3	0	\\
4	0	\\
5	0	\\
6	0	\\
7	0	\\
8	0	\\
9	0	\\
10	0	\\
11	0	\\
12	0	\\
13	0	\\
14	0	\\
15	0	\\
16	0	\\
17	0	\\
18	0	\\
19	0	\\
20	0	\\
};
\end{axis}

\begin{axis}[%
width=0.45\figurewidth,
height=0.3\figureheight,
at={(0.5\figurewidth,0.33\figureheight)},
title = {jump-1024 problem, 6 levels},
xlabel={V-cycle iter.},
ylabel = $\| \vec{x} - \vec{x}^{(n)}\|_{\A}$,
xmin = 0,
xmax = 15,
ymax = 1,
ytick = {1,1.0E-4,1.0E-11,1.0E-13},
xtick = {0,2,9,15,20,30},
yminorticks=true,
ymode = log,
axis background/.style={fill=white}
]
\addplot[
color=mycolor1,
style = solid,
line width=1.0pt,
mark options={solid},
forget plot
]
table[row sep=crcr]{%
0	0.066698707	\\
1	0.000700793	\\
2	3.51E-05	\\
3	2.84E-06	\\
4	2.80E-07	\\
5	3.12E-08	\\
6	3.76E-09	\\
7	4.75E-10	\\
8	6.17E-11	\\
9	8.92E-12	\\
10	3.79E-12	\\
11	3.35E-12	\\
12	3.17E-12	\\
13	3.03E-12	\\
14	2.90E-12	\\
15	2.79E-12	\\
16	2.69E-12	\\
17	2.59E-12	\\
18	2.50E-12	\\
19	2.41E-12	\\
20	2.33E-12	\\
};
\addplot[
color=mycolor2,
style = densely dashed,
line width=1.0pt,
mark options={solid},
forget plot
]
table[row sep=crcr]{%
0	0.066698707	\\
1	0.000701056	\\
2	4.65E-05	\\
3	2.89E-05	\\
4	2.30E-05	\\
5	1.95E-05	\\
6	1.70E-05	\\
7	1.52E-05	\\
8	1.38E-05	\\
9	1.27E-05	\\
10	1.17E-05	\\
11	1.10E-05	\\
12	1.03E-05	\\
13	9.73E-06	\\
14	9.22E-06	\\
15	8.76E-06	\\
16	8.35E-06	\\
17	7.97E-06	\\
18	7.63E-06	\\
19	7.31E-06	\\
20	7.01E-06	\\
};
\addplot[
color=mycolor8,
style = solid,
line width=1.0pt,
mark options={solid},
forget plot
]
table[row sep=crcr]{%
0	0.066698707	\\
1	0.000700793	\\
2	3.51E-05	\\
3	2.84E-06	\\
4	2.80E-07	\\
5	3.12E-08	\\
6	3.76E-09	\\
7	4.75E-10	\\
8	6.16E-11	\\
9	8.14E-12	\\
10	1.09E-12	\\
11	1.60E-13	\\
12	6.51E-14	\\
13	6.19E-14	\\
14	6.16E-14	\\
15	6.15E-14	\\
16	6.17E-14	\\
17	6.16E-14	\\
18	6.18E-14	\\
19	6.14E-14	\\
20	6.17E-14	\\
};

\end{axis}

\begin{axis}[%
width=0.45\figurewidth,
height=0.3\figureheight,
at={(0.5\figurewidth,0\figureheight)},
title = {jump-1024 problem, 6 levels},
xlabel={V-cycle iter.},
ylabel = number of CG iter.,
xtick = {0,2,9,15,20,30},
xmin = 0,
xmax = 15,
yminorticks=true,
axis background/.style={fill=white}
]
\addplot[
color=mycolor1,
style = solid,
line width=1.0pt,
mark options={solid},
forget plot
]
table[row sep=crcr]{%
1	1121	\\
2	985	\\
3	889	\\
4	827	\\
5	733	\\
6	614	\\
7	472	\\
8	238	\\
9	101	\\
10	15	\\
11	0	\\
12	0	\\
13	0	\\
14	0	\\
15	0	\\
16	0	\\
17	0	\\
18	0	\\
19	0	\\
20	0	\\
};
\addplot[
color=mycolor2,
style = densely dashed,
line width=1.0pt,
mark options={solid},
forget plot
]
table[row sep=crcr]{%
1	361	\\
2	303	\\
3	0	\\
4	0	\\
5	0	\\
6	0	\\
7	0	\\
8	0	\\
9	0	\\
10	0	\\
11	0	\\
12	0	\\
13	0	\\
14	0	\\
15	0	\\
16	0	\\
17	0	\\
18	0	\\
19	0	\\
20	0	\\
};
\end{axis}

\end{tikzpicture}%
\caption{
\rv{
Properties of inV-cycle methods with CG as the solver on the coarsest level, which is stopped when the assumption on the absolute coarsest-level accuracy \eqref{eq:thm:eps_def} (approximately) holds with $\epsilon = \theta( 1 - \| \matrx{E} \|_{\A})$, where $\theta=10^{-4}$ (\ref{line:ex2_delta_1e-4}) or  $\theta=10^{-11}$ (\ref{line:ex2_delta_1e-11}). For comparison we also include the $\A$-norm of the error of the exV-cycle method (\ref{line:ex2_backslash}).
}
}
\label{fig:ex2_01_12_23}
\end{figure}

\rv{
The results are summarized in \Cref{fig:ex2_01_12_23}.
After each V-cycle iteration we compute the $\A$-norm of the difference of the exV-cycle and inV-cycle approximations after $n$ V-cycle iterations, i.e., 
\begin{equation}\label{eq:ex2_diff}
\|  \vec{x}^{(n)}_{\ex}  - \vec{x}^{(n)}_{\inex}  \|_{\A}, \quad n=1,2, \ldots, \quad.
\end{equation}
According to estimate \eqref{eq:thm:absolute_aproach_difference_estimate} from \Cref{thm:absolute_approach}, the norm of the difference \eqref{eq:ex2_diff} should be less than 
\begin{equation*}
 \frac{\epsilon \| \matrx{S} \|_{\matrx{A}_0,\A}}{1-\| \matrx{E} \|_{\A}}. 
\end{equation*}
Bounding $\| \matrx{S} \|_{\matrx{A}_0,\A}$ from above by one and considering our choice of $\epsilon$, we get that the difference \eqref{eq:ex2_diff} should be less than $\theta$. Looking at the results, we see that the computed values \eqref{eq:ex2_diff} are slightly less than $\theta$. The estimate \eqref{eq:thm:absolute_aproach_difference_estimate} is accurate for these numerical experiments.
}

\rv{
The convergence of the inV-cycle and exV-cycle methods are approximately the same until they reach the level $\theta$. The $\A$-norm of the error of the inV-cycle method then starts decreasing with a significantly slower rate. At this point the stopping criterion on the coarsest-level is automatically satisfied and the coarsest-level solver is not used. The method perform only smoothing on the fine levels.
}

\rv{
We see that the choice of $\epsilon$, respectively $\theta$, determines the finest-level accuracy of the inV-cycle approximation.
If we look at the number of coarsest-level solver iterations they are decreasing with each V-cycle iteration until they reach zero.
The number of CG iterations performed for the variant with $\theta=10^{-4}$ is significantly smaller than for the variant with $\theta=10^{-11}$.
}

\rv{
The behaviour is analogous for the two problems, the method for the jump-1024 requires significantly more coarsest-level iterations.
}

\subsection{inV-cycle method with absolute coarsest-level stopping criteria}
\label{subsec:ex4_absolute_stopping}
\rv{
In this experiment we study the behaviour of inV-cycle methods with an absolute coarsest-level stopping criteria based on upper bounds of the $\A_0$-norm of the errors.
}

\rv{
We run analogous numerical experiments as in \Cref{subsec:ex2_absolute_ERR}. The only difference is that we stop CG on the coarsest-level using the stopping criterion \eqref{eq:computable_absolute_stopping_crit}, i.e., when
\begin{equation*}
      \eta(\vec{v}_{0,\inex}) \leq \epsilon,
\end{equation*}
where $\eta$ is an upper bound on the $\A_0$-norm of the error of the coarsest-level solver. We again choose $\epsilon =\theta( 1 - \| \matrx{E} \|_{\A} )$, where $\theta=10^{-4}$ or $\theta=10^{-11}$.
We consider two choices of $\eta$. First, the residual based upper bound \eqref{eq:res_upper_bound}. We label this variant as RES. We approximate the term $\| \A_0^{-1}\|$ using the MATLAB function \texttt{eigs}. Second, the Gauss-Radau upper bound on the $\A_0$-norm of the error in CG stated in \cite[second inequality in (3.5) with updating formula for a coefficient (3.3)]{Meurant2023}.
This upper bound is based on the interpretation of CG as a procedure for computing a Gauss-Radau quadrature approximation to a Riemann-Stieltjes integral. To compute this upper bound we need an lower bound on the smallest eigenvalue of the matrix $\A_0$. We approximate the smallest eigenvalue of $\A_0$ using the MATLAB \texttt{eigs} function and use its $1-10^{-3}$ multiple as the lower bound. We label this variant as GR. For comparison we include in the plots the results computed in \Cref{subsec:ex2_absolute_ERR} where CG is stopped on the coarsest-level when inequality \eqref{eq:thm:eps_def} (approximately) holds. We label this variant as ERR.
}

\begin{figure}
\setlength\figureheight{18cm}
\setlength\figurewidth{0.99\textwidth} 
\input{img/inV-cycleCGAbsoluteCoarsestLevelStoppingCriteria.tikz}
\caption{
\rv{
Properties of inV-cycle methods with CG as the solver on the coarsest level, which is stopped by an absolute criterion based on upper bounds of the $\A_0$-norm of the errors; variant ERR with $\theta=10^{-4}$ (\ref{line:ex4_ERR_1e-4}) or $\theta=10^{-11}$ (\ref{line:ex4_ERR_1e-11}), variant GR with $\theta=10^{-4}$ (\ref{line:ex4_GR_1e-4}) or  $\theta=10^{-11}$ (\ref{line:ex4_GR_1e-11}), variant RES with $\theta=10^{-4}$ (\ref{line:ex4_RES_1e-4}) or  $\theta=10^{-11}$ (\ref{line:ex4_RES_1e-11}). For comparison we also include the $\A$-norm of the error and Euclidean norm of residual of the exV-cycle method (\ref{line:backslash_ex4}).
}
}
\label{fig:ex4}
\end{figure}
\rv{
We run the V-cycle method starting with a zero initial approximate solution and stop after $15$ V-cycle iterations.
The results are summarized in \Cref{fig:ex4}. After each V-cycle iteration we compute the $\A$-norm of the difference of the exV-cycle and inV-cycle approximations after $n$ V-cycle iterations \eqref{eq:ex2_diff}. According to the discussion in \Cref{sec:new_stopping_criteria} and the choice of $\epsilon$, the norm of the difference \eqref{eq:ex2_diff} should be less than~$\theta$. Looking at the results we see that all values \eqref{eq:ex2_diff} are lower than the corresponding~$\theta$. 
We see that estimate \eqref{eq:rel_diff_abs_stop_E_norm} is the most accurate for the variant ERR and the loosest for the variants RES.
When performing the experiments we observed that the Gauss-Radau upper bound on the $\A_0$-norm of the error used in the GR variants is more accurate than the residual based estimate \eqref{eq:res_upper_bound} used in the RES variants.
The more accurate the upper bound on the $\A_0$-norm of the error on the coarsest-level is used in the stopping criterion the more accurate estimate \eqref{eq:rel_diff_abs_stop_E_norm} is and the less CG iterations on the coarsest-level are performed.
}

\rv{
Looking at the $\A$-norms of the error, we see that the variants GR and RES with stopping criteria based on the upper bounds of the $\A_0$-norm of the coarsest-level errors have analogous convergence behavior as the variant ERR with stopping criteria based on the $\A_0$-norm of the coarsest-level errors. 
}

\rv{
Based on these experiments, we believe that automatic satisfaction of the coarsest-level criteria can be used as a heuristic indicator that the $\A$-norm of the error on the finest level is at the level of $\theta$. Another heuristic indicator that we reached the desired finest-level accuracy might be a stagnation of the norm of the finest-level residual.
}

\subsection{Performance of inV-cycle methods with absolute coarsest-level stopping criteria}
\rv{
In this experiment, we evaluate the performance of inV-cycle methods with an absolute coarsest-level stopping criteria considered in \Cref{subsec:ex4_absolute_stopping}. 
}

\rv{
We consider the same problems and analogous V-cycle methods. The only difference is that we don't use a computed approximation of $ \| \matrx{E} \|_{\A}$ but assume that $ \| \matrx{E} \|_{\A} < 2/3$ for both problems. The assumption $\| \matrx{E} \|_{\A}<2/3$ should be a valid assumption for most of the well set up V-cycle methods. For difficult problems it may be safer to consider it closer to one.
Our goal is to compute approximations whose $\A$-norm of the error is approximately at the level of $10^{-4}$ and $10^{-11}$, respectively.
According to the discussion in \Cref{sec:new_stopping_criteria} we choose $\epsilon = (1-2/3)\theta$, where $\theta = 10^{-4}$ and $\theta = 10^{-11}$. 
}

\begin{figure}
\centering
\includegraphics[width=0.9\textwidth]{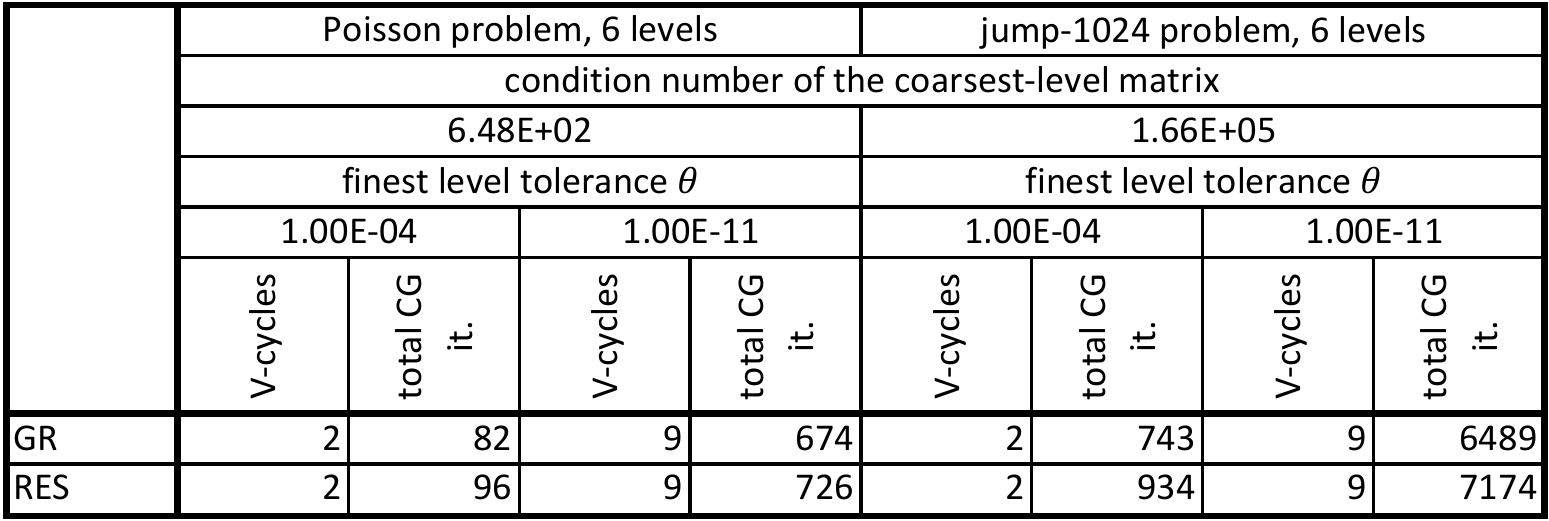}
\caption{
\rv{
Properties of inV-cycle methods with CG as the solver on the coarsest level, which is stopped by the absolute criteria based on upper bounds of the $\A_0$-norm of the errors.
}
}
\label{fig:ex5_motiv}
\end{figure}

\rv{
We run the V-cycle method starting with a zero initial approximate solution and stop when the $\A$-norm of the error is (approximately) lower than $10^{-4}$ and $10^{-11}$ for the variants with $\theta=10^{-4}$ and $\theta=10^{-11}$, respectively.
For both problems the exV-cycle method requires $2$ and $9$ V-cycle iterations to reach the desired finest-level accuracy $10^{-4}$ and $10^{-11}$, respectively. The results of the inV-cycle methods are summarized in \Cref{fig:ex5_motiv}.
}

\rv{
We see that the inV-cycle methods converge to the desired accuracy in the same number of V-cycle iterations as the exV-cycle methods. The goal of the coarsest-level stopping strategy is thus satisfied. The methods works well for both problems with the same choice of the parameter $\epsilon$. The variants RES, require more CG iterations on the coarsest level than the variants GR.
}

\rv{
We may compare the total number of CG iterations in the variants GR and RES with the total number of CG iterations in the variants with a relative residual stopping criterion in \Cref{fig:motiv}. We see that the number of total CG iterations in the GR and RES variants are not the lowest possible, such that an inV-cycle method converges to the desired accuracy in the same number of V-cycles as the exV-cycle method, but they also aren't substantially high.
}

\begin{figure}
\centering
\includegraphics[width=0.9\textwidth]{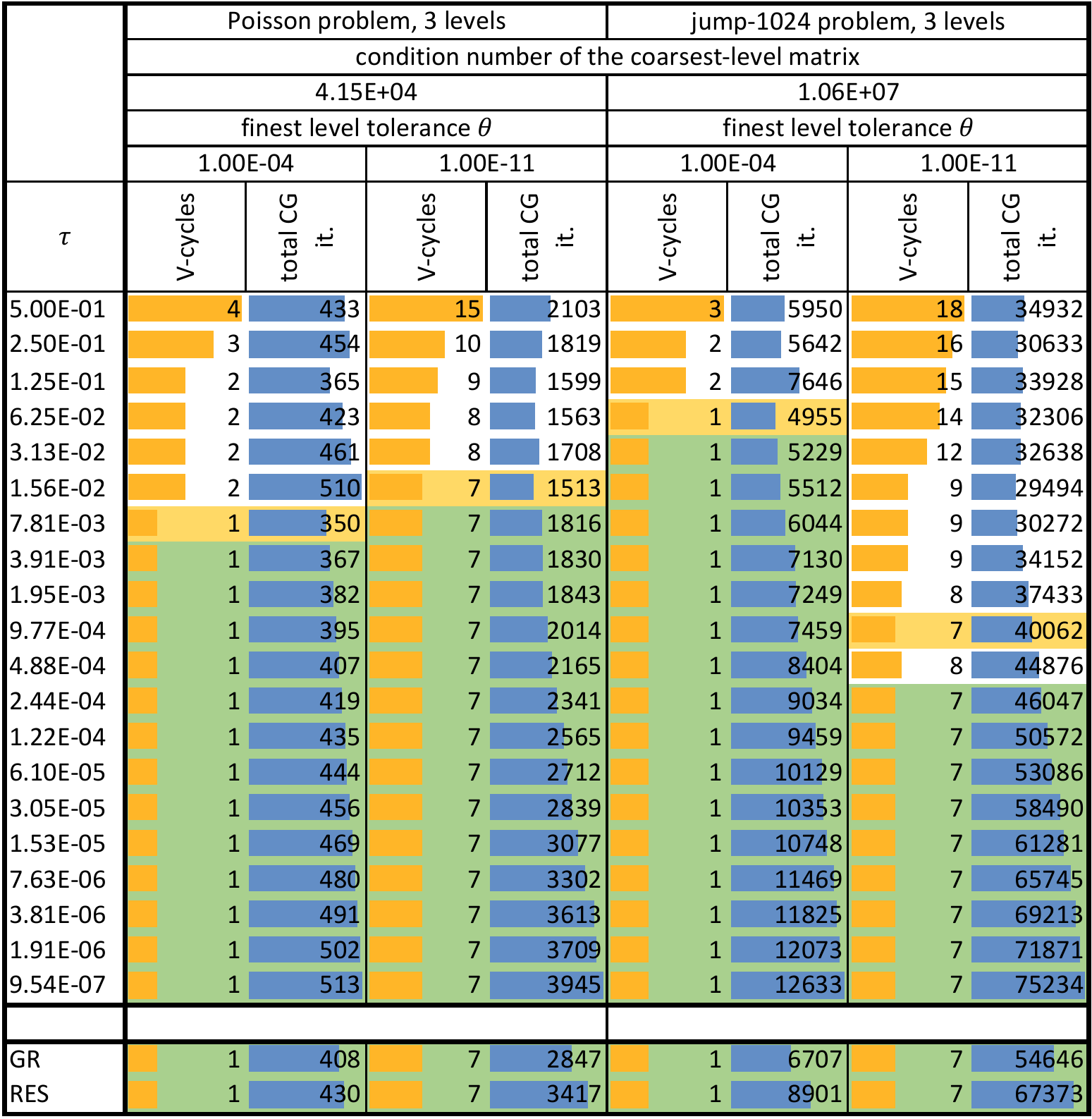}
\caption{
\rv{
Properties of inV-cycle methods with CG as the solver on the coarsest level, which is stopped by a relative residual criterion with various tolerance $\tau$, or by an absolute criterion based on upper bounds of the $\A_0$-norm of the errors; variants GR and RES.
The bright yellow and green color highlight variants that converge in the same number of V-cycles as the exV-cycle method. The bright yellow variants achieve this in the least total number of CG iterations on the coarsest-level.
}
}
\label{fig:ex5_3levels}
\end{figure}

\rv{
To see how the coarsest-level stopping strategy may be affected by the change of the size of the coarsest-level problem and the change of the number of levels in the V-cycle method we run experiments where we consider the same problem on the finest level, but just three level V-cycle methods. The size of the coarsest-level problems is $101761$ DoFs. The results are summarized in \Cref{fig:ex5_3levels}.
}

\rv{
We see analogous behavior as in the experiment with six level V-cycle methods. The variants GR and RES converge to the desired accuracy in the same number of V-cycle iterations as the exV-cycle methods.
}

\rv{
The main benefit of the stopping strategy is that we don't have to try different parameters for different problems or when we want to reach different finest-level tolerances or when the size of the coarsest-level problem changes. The parameter $\theta$ is chosen the same as the finest-level tolerance we are aiming for.
}

\section{Conclusions and open problems}
\label{sec:conclusion}
\rv{
In this paper we present an approach to analyzing the effects of approximate coarsest-level solves on the convergence of the V-cycle method for SPD problems.
We use the results to give an answer to the question of how the choice of tolerance in the relative residual stopping criterion for the coarsest-level solver may affect the convergence of the V-cycle method.
We present novel coarsest-level stopping criterion which we may use to control the difference between the computed approximation and the approximation which would be computed by the exV-cycle method.
This coarsest-level stopping criterion may thus be set up such that the method converges to a chosen finest-level accuracy in (nearly) the same number of V-cycle iterations as the exV-cycle method.
The stopping strategy achieves this goal in various numerical experiments. In a future work we would like to test this coarsest-level stopping strategy within the algebraic multigrid methods.
}

\rv{
In this work we focus on the use of multigrid methods as a standalone solver. Multigrid methods are, however, also frequently used as a preconditioner for a Krylov subspace method.
It would be interesting to investigate how the results obtained in this paper could be utilized in this setting. 
In general an inV-cycle method would have to be applied as a flexible preconditioner.
}

\rv{
Other open problems include the generalization to non-symmetric problems or to other multigrid schemes such as the W-cycle scheme or the full multigrid scheme.
}

\appendix
\section{Numerical approximation of $\| \matrx{E} \|_{\A}$}
\label{appendix}
\rv{
In this section we describe a procedure for numerical approximation of the $\A$-norm of the error propagation matrix $\matrx{E}$ of the exV-cycle scheme. We consider an exV-cycle scheme where the pre- and post- smoothing is each accomplished by one iteration of the symmetric Gauss-Seidel method. Thanks to the use of the symmetric Gauss-Seidel smoother the matrix $\matrx{E}$ is symmetric and there exist a symmetric matrix $\matrx{B}$ such that  $\matrx{E} = \matrx{I} - \matrx{B}^{-1}\matrx{A}$; see, e.g., \cite{Xu1992}.
Then
\begin{equation*}
\| \matrx{E} \|_{\A} =\| \matrx{I} - \matrx{B}^{-1}\matrx{A} \|_{\A} =  \| \matrx{A}^{\frac{1}{2}}(\matrx{I} -  \matrx{B}^{-1} \matrx{A} )\matrx{A}^{-\frac{1}{2}} \| = \| \matrx{I} - \matrx{A}^{\frac{1}{2}} \matrx{B}^{-1} \matrx{A}^{\frac{1}{2}} \|.    
\end{equation*}
Since the matrices $\matrx{A}^{\frac{1}{2}} \matrx{B}^{-1} \matrx{A}^{\frac{1}{2}}$ and $\matrx{B}^{-1} \matrx{A}$ have the same eigenvalues there holds
\begin{equation*}
 \| \matrx{I} - \matrx{A}^{\frac{1}{2}} \matrx{B}^{-1} \matrx{A}^{\frac{1}{2}} \| = \| \matrx{I} -  \matrx{B}^{-1} \matrx{A} \|,    
\end{equation*}
and consequently $\| \matrx{E} \|_{\A} =\| \matrx{E} \|$. We compute it using MATLAB function \texttt{eigs} (with the largest eigenvalue option) applied to the function
\begin{equation*}
     \vec{x} \mapsto \vec{x} - \mathbf{V}( \matrx{A}_{0:J}, \matrx{M}_{1:J}, \matrx{N}_{1:J}, \matrx{P}_{1:J}, \matrx{A}\vec{x},\vec{0},J).
\end{equation*}
}

\section{\rev{Derivation of inequalities \eqref{eq:error_A_norm_by_residual_Eucl_norm_bound_upper} and  \eqref{eq:error_A_norm_by_residual_Eucl_norm_bound_lower}}}\label{apendix:inequalities}
\rev{In this section we present derivations of inequalities \eqref{eq:error_A_norm_by_residual_Eucl_norm_bound_upper} and  \eqref{eq:error_A_norm_by_residual_Eucl_norm_bound_lower} used in \Cref{sec:effects_rel_res}, i.e.,
\begin{align*}
\|\matrx{A}_0^{-1} \|^{-\frac{1}{2}} \|\vec{v}_0 - \vec{v}_{0,\inex}  \|_{\matrx{A}_0} 
&\leq   \| \vec{f}_0 - \matrx{A}_0 \vec{v}_{0,\inex} \|,
\\ \|  \vec{b} -  \A\vec{x}^{\prev} \| 
  & \leq   \| \matrx{A}\|^{\frac{1}{2}}   \| \vec{x} - \vec{x}^{\prev} \|_{\matrx{A}}.
\end{align*}
Using that $\A_0\vec{v}_0 = \vec{f}_0$ and that $\A_0$ is SPD we have
\begin{align*}
\|\vec{v}_0 - \vec{v}_{0,\inex}  \|^2_{\matrx{A}_0} 
& = (\vec{v}_0 - \vec{v}_{0,\inex})^{\top} \matrx{A}_0 (\vec{v}_0 - \vec{v}_{0,\inex}) 
\\& = (\matrx{A}^{-1}_0 (\vec{f}_0 - \matrx{A}_0 \vec{v}_{0,\inex}))^{\top} \matrx{A}_0 (\matrx{A}^{-1}_0(\vec{f}_0 - \matrx{A}_0\vec{v}_{0,\inex} )) 
\\& = ( \vec{f}_0 - \matrx{A}_0 \vec{v}_{0,\inex})^{\top} \matrx{A}^{-1}_0\matrx{A}_0 \matrx{A}^{-1}_0 (\vec{f}_0 - \matrx{A}_0\vec{v}_{0,\inex} )  
\\& = 
( \vec{f}_0 - \matrx{A}_0 \vec{v}_{0,\inex})^{\top} \matrx{A}^{-1}_0 (\vec{f}_0 - \matrx{A}_0\vec{v}_{0,\inex} )
 \leq \| \matrx{A}^{-1}_0 \|   \| \vec{f}_0 - \matrx{A}_0 \vec{v}_{0,\inex} \|^2 ,
\end{align*}
which yields the first inequality.
The second inequality can be derived using that $\A \vec{x} = \vec{b}$ and that $\A$ is SPD
\begin{align*}
\|  \vec{b} -  \A\vec{x}^{\prev} \|^2 
&= (\vec{b} -  \A\vec{x}^{\prev})^{\top} (\vec{b} -  \A\vec{x}^{\prev})
=(\A(\vec{x} -  \vec{x}^{\prev}))^{\top} ((\A(\vec{x} -  \vec{x}^{\prev})) 
\\&= (\vec{x} -  \vec{x}^{\prev})^{\top} \A^{\frac{1}{2}} \A \A^{\frac{1}{2}} (\vec{x} -  \vec{x}^{\prev})
\\& \leq \| \A \| (\vec{x} -  \vec{x}^{\prev})^{\top} \A(\vec{x} -  \vec{x}^{\prev}) = \| \A \| \| \vec{x} -  \vec{x}^{\prev}\|^2_{\A}.
\end{align*}}

\section*{Acknowledgments}
The authors wish to thank Petr Tichý for his useful comments on error estimation in CG and Jaroslav Hron for his suggestions when generating the system matrices in FEniCS. The authors acknowledge the support of the Erasmus+ program that enabled Petr Vacek to spend the Winter semester 2021-2022 at Trinity College Dublin. During this visit the basis of the paper was developed.

\bibliographystyle{siamplain}
\bibliography{references}
\end{document}